\documentclass[letterpaper,twoside,11pt]{article}
\usepackage{authblk,amsmath}

\pdfoutput=1

\usepackage[letterpaper,left=1in,right=1in,top=1in,bottom=1in]{geometry}

\newcommand{\myauthor}{Alain Connes and  Caterina Consani}

\newcommand{\mytitle}{On the metaphysics of F$_1$}

\author{\myauthor\footnote{Partially supported by the Simons Foundation collaboration grant n. 691493}}
\title{\mytitle}

\date{\small To Yuri Ivanovich Manin, in memory.}

\usepackage{mathrsfs}
\usepackage{graphicx}
\usepackage{rotating}
\usepackage{float}
\usepackage{caption}
\usepackage{pdflscape}
\usepackage{url}
\usepackage[all,arc,2cell]{xy}
\UseAllTwocells
\usepackage{enumerate}
 \usepackage{lineno}
\usepackage{stix}
\usepackage{accents}
\usepackage{makecell}
\usepackage{relsize}
\usepackage{tikz}
\usepackage{amsmath,amsthm}
\usepackage{caption}
\usepackage{subcaption}

\usepackage{microtype}
\usepackage{amsthm}
\usepackage{stmaryrd}

\usepackage{marginnote}

\usepackage{spectralsequences}
\usepackage[nomessages]{fp}
\usepackage{ifthen}

\usepackage[pdfstartview=FitH, pdfauthor={\myauthor},
            pdftitle={\mytitle},
            colorlinks,
            linkcolor=reference,
            citecolor=citation,
            urlcolor=e-mail,
            backref]{hyperref}
\usepackage[nameinlink, capitalise]{cleveref}

\usepackage{fancyhdr}
\usepackage{color}
\definecolor{todo}{rgb}{1,0,0}
\definecolor{conditional}{rgb}{0,1,0}
\definecolor{e-mail}{rgb}{0,.40,.80}
\definecolor{reference}{rgb}{.20,.60,.22}
\definecolor{mrnumber}{rgb}{.80,.40,0}
\definecolor{citation}{rgb}{0,.40,.80}

\let\oldmarginpar\marginpar
\renewcommand\marginpar[1]{\-\oldmarginpar[\raggedleft\footnotesize #1]%
{\raggedright\footnotesize #1}}

\def\A{{\mathbb A}}
\def\B{{\mathbb B}}
\def\F{{\mathbb F}}
\def\N{{\mathbb N}}
\def\Q{{\mathbb Q}}
\def\R{{\mathbb R}}
\def\Z{{\mathbb Z}}
\def\C{{\mathbb C}}

\def\cA{{\mathcal A}}

\def\cP{{\mathcal P}}
\def\cO{{\mathcal O}}

\def\Spec{{\rm Spec\,}}
\newcommand{\sss}{{{\mathbb S}}}
\def\spm{{\sss[\pm 1]}}
\def\Se{\mathfrak{ Sets}}
\def\Ses{{\Se_*}}
\def\Hom {{\rm{Hom}}}

\def\gop{{\Gamma^{o}}}

\def\spzb{{\overline{\Spec\Z}}}

\def\zmax{{\Z_{\rm max}}}
\def\rmax{\R_+^{\rm max}}
\def\gop{{\Gamma^{\rm op}}}
\def\spzb{{\overline{\Spec\Z}}}
\def\wnt{{\widehat{\N^{\times}}}}
\def\aarith{{\mathscr A}}
\def\arith{{(\wnt,\cO)}}
\def\fr{{\rm Fr}}
\def\hatz{{\hat\Z^\times}}
\def\cdim{{{\mbox{Dim}_\R}}}
\def\vsp{\vspace{.05in}}
\def\nt{\N^{\times}}
\def\tdim{{{\mbox{dim}_{\rm top}}}}
\def\cE{{\mathcal E}}
\def\div{{\rm Div}}
\def\Aut{{\rm Aut}}
\def\gam{{\Gamma\Ses}}
\def\catmo{{\mathfrak{Mod}}}

\newcommand{\ie}{{\it i.e.\/}\ }

\newcommand{\opcit}{{\it op.cit.\/}\ }
\newcommand{\cf}{{\it cf.}}

\newcommand{\no}{\noindent}

\numberwithin{equation}{section}
\theoremstyle{plain}
\newtheorem{theorem}[equation]{Theorem}
\newtheorem*{theorem*}{Theorem}
\newtheorem{lemma}[equation]{Lemma}

\newtheorem{proposition}[equation]{Proposition}

\newtheorem{corollary}[equation]{Corollary}

\newtheoremstyle{named}{}{}{\itshape}{}{\bfseries}{.}{.5em}{#1 \thmnote{#3}}
\theoremstyle{named}

\theoremstyle{definition}
\newtheorem{definition}[equation]{Definition}

\newtheorem{example}[equation]{Example}

\newtheorem{remark}[equation]{Remark}

\begin{document}

\maketitle


\begin{abstract}
    \noindent 
    In the present paper, dedicated to Yuri Manin, we investigate the general notion of rings of $\sss[\mu_{n,+}]$--polynomials and relate this concept to the known notion of number systems. The Riemann-Roch theorem for the ring $\Z$ of the integers that we obtained recently uses the understanding of $\Z$ as a ring of polynomials $\sss[X]$ in one variable over the absolute base $\sss$, where $1+1=X+X^2$. The absolute base $\sss$ (the categorical version of the sphere spectrum) thus turns out to be a strong candidate for the incarnation of the mysterious $\F_1$.


    \paragraph{Mathematics Subject Classification 2010:}
 \paragraph{Mathematics Subject Classification 2010:}
\href{http://www.ams.org/mathscinet/msc/msc2010.html?t=14Fxx&btn=Current}{14C40},
    \href{http://www.ams.org/mathscinet/msc/msc2010.html?t=14Hxx&btn=Current}{14G40},
    \href{http://www.ams.org/mathscinet/msc/msc2010.html?t=14Kxx&btn=Current}{14H05},
    \href{http://www.ams.org/mathscinet/msc/msc2010.html?t=11Rxx&btn=Current}{11R56},
    \href{http://www.ams.org/mathscinet/msc/msc2010.html?t=13Fxx&btn=Current}{13F35},
    \href{http://www.ams.org/mathscinet/msc/msc2010.html?t=18Nxx&btn=Current}{18G60},
    \href{http://www.ams.org/mathscinet/msc/msc2010.html?t=19Dxx&btn=Current}{19D55}.
\end{abstract}


\section{Introduction}
\label{intro}
\begin{quote}
Les math\' ematiciens du xvi-\`eme si\`ecle avaient coutume de parler de la ``m\' etaphysique du calcul infinit\'esimal", de la ``m\' etaphysique de la th\' eorie des \' equations". Ils entendaient par l\`a un ensemble d'analogies vagues, difficilement saisissables et difficilement formulables, qui n\' eanmoins leur semblaient jouer un r\^ole important \`a un moment donn\' e dans la recherche et la d\' ecouverte math\' ematiques.\hspace{.95in} (A. Weil, De la m\' etaphysique aux math\' ematiques, 1960, \cite{Weil1})
\end{quote}
Yuri Manin, to whose memory we dedicate this article, first recognized in \cite{Man-zetas} the importance of developing a theory of  ``absolute coefficients'' in arithmetic geometry,   independently of the early ideas proposed by R. Steinberg \cite{Steinberg} and J. Tits \cite{Tits} in the context of Chevalley groups. In arithmetic, for number fields, the goal is to provide the geometric counterpart to the construction that A. Weil used in his proof of the Riemann hypothesis for function fields. The search for a close analogy between number fields and function fields of curves in positive characteristic induced Manin to postulate the existence of the absolute point ``$\Spec \F_1$,'' over which one could apply Weil's strategy to the study of the Riemann zeta function. For the algebraic scheme $\Spec \Z$, one would then use the spectrum of the tensor product ``$\Z \otimes_{\F_1} \Z$'' as a substitute for the self-product of a curve over (the spectrum of) a finite field. \newline 
Manin always advocated the fruitfulness of unexpected interactions between different approaches to a mathematical problem. In Sections \ref{cincarnation} and \ref{sectcoeffs} we shall discuss two of such unexpected occurrences, in fact two pillars of our joint work in the past fifteen years. Section \ref{cincarnation} is about the hypothetical curve\footnote{we reserve throughout the symbol $\bf C$ for this entity} $\bf C$ that we propose as the absolute geometric entity. Section \ref{sectcoeffs} concerns instead the absolute coefficients. The aim of this paper is to sponsor $\sss$  the most basic combinatorial form of the sphere spectrum and  of an $\sss$-algebra,  as the most natural candidate for the absolute coefficients (aka $\F_1$). We claim that this algebra is the absolute ``field'' of constants over which $\Z$ becomes a ring of polynomials in one variable. 
 This point of view is supported by the Riemann-Roch theorem for the ring $\Z$  recently proved in \cite{RR1}, whose formula shows that the genus of $\spzb$ is  zero. In an earlier result on the same topic  \cite{RR}, the integers were considered as polynomials over $\spm$ with generator $X=3$. This fact is based on the balanced ternary numeral system\footnote{An early occurrence of this numeral system is found in the 1544 book ``Arithmetica integra" of Michael Stifel.} which provides a balanced signed-digit representation of the integers as finite sums of powers of the ``variable''  $X=3$ with coefficients in the set $\{-1,0,1\}$ underlying the pointed multiplicative monoid $\mu_{2,+}$ of quadratic roots of unity. The new version of the Riemann-Roch theorem for the ring $\Z$  in \cite{RR1} simplifies the earlier version \cite{RR} and it also reconciles the formula (and our understanding of this subject) with the classical number theoretic viewpoint. Indeed,   in the analogy between number fields and curves over finite fields, the field $\Q$  has genus zero \cite{Weil39} and it is singled out as the only field   contained in any other number field. The view of $\Z$  as a ring of polynomials over the absolute base $\sss$ selects the generator $X=-2$. The key fact  is that any integer can be uniquely written as a sum of powers of $-2$ \cite{Knuth}.

\no The above special cases of generators $X$ for rings over finite spherical $\sss$-algebras justify a systematic and broader study of rings of $\sss[\mu_{n,+}]$--polynomials. 
  In Section \ref{sec5} we introduce the general notion of rings of $\sss[\mu_{n,+}]$--polynomials in one and several variables. Let $n>0$ be an integer, $\mu_n$ the multiplicative group of $n$-th roots of $1$ and  $\sss[\mu_{n,+}]$ the spherical $\sss$-algebra of the (pointed) monoid $\mu_{n,+}=\mu_n\cup \{0\}$. We recall that morphisms of $\sss$-algebras $\sss[\mu_{n,+}]\to HR$ ($R$ being a ring) correspond bijectively to group homomorphisms $\iota: \mu_n\to R^\times$ \cite{CCprel}. Let $\cP(\mu_n)$ be the subset of the set $\left(\mu_n\cup \{0\}\right)^\N$ of sequences with only finitely many non-zero terms. By definition, an element $X\in R$ is an  $\sss[\mu_{n,+}]$-generator if and only if the evaluation map $\sigma:\cP(\mu_n)\to R$, $\sigma((\alpha_j))=\sum_j \iota(\alpha_j) \, X^j$ is bijective. Proposition \ref{sprop2} shows that the pair $(R, X)$ of a ring of $\sss[\mu_{n,+}]$--polynomials in one variable is uniquely specified, up to isomorphism, by  the map $h:\mu_n\to \cP(\mu_n)$,  which, in turn, is uniquely defined by the equality $\sigma(h(\xi))=\iota(\xi)+1$. In Section \ref{sec6}  we give several examples of rings  of $\sss[\mu_{n,+}]$--polynomials based on some known number systems. We refer to \cite{BBLT} for a survey on systems of numerations and for references therein contained, but we claim no exhaustiveness. Conceptually, the examples of rings of $\sss[\mu_{n,+}]$-polynomials discussed in this article provide an explicit bridge between the $p$-adic and the complex world. At the geometric level, the rings of polynomials are naturally related to the projective line $\mathbb P^1$, and the evaluation at the points $0$ and $\infty$ of    
  $\mathbb P^1$ yields, after completion, the following refinement (the lower line) of a classical diagram (upper line). In the upper line, $K$ is the field of fractions of the $p$-typical Witt ring of the algebraic closure of $\F_q$ ($q=p^\ell$) and $\overline K$ is its algebraic closure.
\[
\renewcommand{\arraystretch}{1.3}
\begin{array}[c]{ccccccccc}
\overline\F_q&\stackrel{\pi}{\twoheadleftarrow}&W(\overline\F_q)&\hookrightarrow& \overline K&\supset&\overline\Q&\subset&\C\\
\rotatebox{90}{$\subset$}&&\rotatebox{90}{$\subset$}&&\rotatebox{90}{$\subset$}&&\rotatebox{90}{$\subset$}&&\rotatebox{90}{$=$}\\
\F_q&\stackrel{\pi}{\twoheadleftarrow}&W(\F_q)&\hookrightarrow&W(\F_q)[\eta]&\hookleftarrow& R[X^{-1}]&\hookrightarrow& \C
\end{array}
 \]
\no In the lower line, $X$ is a $\sss[\mu_{n,+}]$-generator of the ring $R$ where $n+1=q$.  $R[X^{-1}]$ is the ring of Laurent polynomials; the map to $\C$ is  the inclusion of $R[X^{-1}]$ in $\C$ by specialization of $X$, obtained by solving  the equations $\sigma(h(\xi))=\iota(\xi)+1, \, \xi \in \mu_n$, and using  the canonical embedding  $\mu_{n,+}\subset \C$. The map from $R[X^{-1}]$ to the finite extension $W(\F_q)[\eta]$ is obtained from the canonical inclusion of $R$ in the projective limit $\varprojlim R_n$ (see Proposition \ref{sprop2}).

\no The general theory of rings of $\sss[\mu_{n,+}]$-polynomials, together with the role of the absolute base $\sss$ in the formulation of the Riemann-Roch theorem  \cite{RR1}, suggest the following refinement of the definition of the Arithmetic Site. Originally, this space was defined by the pair of the arithmetic topos $\wnt$  and the structure sheaf given by the Frobenius action of $\nt$ on the tropical semiring $\zmax$ \cite{CCas}. The role of the field of constants is here played by the Boolean semifield $\B$. The   development of this paper evidently hints to a replacement of the structure sheaf $\zmax$ by the sheaf of $\sss$-algebras obtained from the Frobenius action $X\mapsto X^n$ of $\nt$ on the spherical algebra $\sss[X]$. This new version of $\sss$-arithmetic site provides simultaneously a natural base both at the coefficients  and at the geometric levels.  The topos $\wnt$ is the geometric incarnation of the $\lambda$-operations in the theory of $\lambda$-rings \cite{Borgers} in the context of geometry over $\F_1$. We expect that throughout a suitable understanding of the ``algebraic closure" $\overline \F_1$ of the absolute coefficients one may relate the space of points of the $\sss$-arithmetic site over $\overline \F_1$  with the (points of the) curve $\bf C$ whose structure is recalled in Section \ref{cincarnation}.

Finally, these results also point out to the open and interesting question of the classification of rings of $\sss[\mu_{n,+}]$--polynomials in several variables which  pursues  the intuitive statement of Yuri Manin  \cite{Man-zetas}:
  \begin{quote}
      \emph{The central question we address can be provocatively put as follows: if numbers are similar to polynomials in one variable over a finite field, what is the analog of polynomials in several variables? Or, in more geometric terms, does there exist a category in which one can define ``absolute Descartes powers'' $\Spec\Z\times\cdots\times \Spec\Z$?}
  \end{quote}

  \section{Adelic and topos theoretic incarnation of $\bf C$} \label{cincarnation}
A first connection  between  Manin's point of view on  $\F_1$ and a seemingly unrelated topic takes place  as a by-product of the relations between C. Soul\'e  perspective on varieties over $\F_1$ (named ``Critical Realism'' in \cite{M}) -- motivated by Manin \cite{Man-zetas} (\cf~\S 1.5) --
and the work of the first author \cite{Co-zeta} on the trace formula in noncommutative
geometry and the zeros of the Riemann zeta function. In \cite{Soule}, Soul\'e introduced the following zeta function of a variety $X$ over $\F_1$  
\begin{equation}\label{zetadefn}
\zeta_X(s):=\lim_{q\to 1}Z(X,q^{-s}) (q-1)^{N(1)},\qquad s\in\R
\end{equation}
using the {\it polynomial} counting function $N(x)\in \Z[x]$ associated with $X$ and  the Hasse-Weil exponential series
\begin{equation}\label{zetadefn1}
Z(X,T) := \exp\left(\sum_{r\ge 1}N(q^r)\frac{T^r}{r}\right).
\end{equation}
All the examples of varieties considered in \opcit are rational. Thus,  the existence of an underlying curve $\bf C$ related, in a similar manner, to the Riemann zeta function is subordinated to finding a function $N(q)$ (highly non-polynomial!)  that produces, through the use of \eqref{zetadefn},  the complete Riemann zeta function  $\zeta_\Q(s)=\pi^{-s/2}\Gamma(s/2)\zeta(s)$. This is a non-trivial problem since  classically, $N(1)$ in the above formula inputs the Euler characteristic of the geometric space. Thus one might be induced to expect\footnote{the number of zeros of $\zeta_\Q$ is infinite, and so is the dimension of the (mysterious) cohomology $H^1(\bf C)$} that since for the Riemann zeta-function one ought to have $N(1)=-\infty$,   the use of \eqref{zetadefn}  should be precluded, and with it also the expectation that $N(q)\geq 0$ for $q\in (1,\infty)$. There is, in fact,   a  natural way to by-pass this problem by  applying  the logarithmic derivatives to both sides of  \eqref{zetadefn} and then observing that the right-hand side determines the Riemann sums of an integral \cite{CC0, CC1}. In this way, in place of  \eqref{zetadefn} one considers the equation:
$
\frac{\partial_s\zeta_N(s)}{\zeta_N(s)}=-\int_1^\infty  N(u)\, u^{-s}d^*u,
$
where $d^*u:=du/u$.
This integral formula  produces the following one  for the sought for counting function $N(q)$ associated with $\bf C$:
\begin{equation}\label{special}
   \frac{\partial_s\zeta_\Q(s)}{\zeta_\Q(s)}=-\int_1^\infty  N(u)\, u^{-s}d^*u.
\end{equation}
The above equation admits a meaningful solution expressable in terms of the {\it distribution} 
\begin{equation}\label{Nu}
    N(u)=\frac{d}{du}\varphi(u)+ \kappa(u), \qquad \varphi(u):=\sum_{n<u}n\,\Lambda(n),
\end{equation}
where  $\kappa(u)$ is the distribution that appears in the Riemann-Weil explicit formula
$$
\int_1^\infty\kappa(u)f(u)d^*u=\int_1^\infty\frac{u^2f(u)-f(1)}{u^2-1}d^*u+cf(1)\,, \qquad c=\frac12(\log\pi+\gamma).
$$
One shows that the distribution $N(u)$ is positive on $(1,\infty)$, and when written in terms of the non-trivial zeros $\rho\in Z$ of the Riemann zeta function,  it is given, in complete analogy with its counterpart holding in the function field case, by
\begin{equation}\label{fin2}
    N(u)=u-\frac{d}{du}\left(\sum_{\rho\in Z}{\rm order}(\rho)\frac{u^{\rho+1}}{\rho+1}\right)+1,
\end{equation}
where the derivative is taken in the sense of distributions. The value at $u=1$ of the  term
 $\displaystyle{\omega(u)=\sum_{\rho\in Z}{\rm order}(\rho)\frac{u^{\rho+1}}{\rho+1}}$ is given  by
$\frac 12+ \frac \gamma 2+\frac{\log4\pi}{2}-\frac{\zeta'(-1)}{\zeta(-1)}
$. 
 \begin{figure}[H]
\begin{center}
\includegraphics[scale=1]{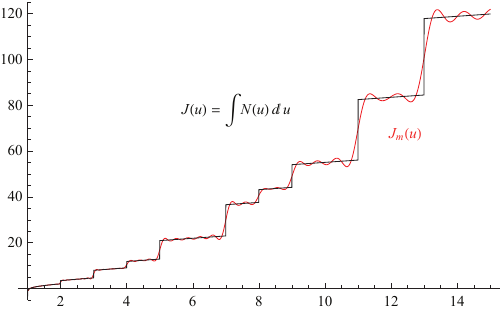}
\end{center}
\caption{Graph of a primitive $J(u)$  of the counting distribution $N(u)$. One has $J(u)\to -\infty$ when $u\to 1$. The wiggly graph is the approximation of $J(u)$ obtained using the symmetric set $Z_m$ of the first $2m$ zeros to perform the sum $J_m(u)=\frac{u^2}{2}-\sum_{Z_m}{\rm order}(\rho)\frac{u^{\rho+1}}{\rho+1}+u$.}
\label{figcounting}     
\end{figure}
The tension between the positivity of the distribution $N(q)$ for $q>1$, and the expectation that its value $N(1)$ ought to be $N(1)=-\infty$ is  resolved by implementing the theory of distributions. Indeed, even though $N(u)$ is {\it finite} as a distribution,  when one looks at it as a function, its value at $q=1$ is formally given by
$$
N(1)=2-\lim_{\epsilon\to 0}\frac{\omega(1+\epsilon)-\omega(1)}{\epsilon}\sim-\frac 12 E \log E,\qquad \ E=\frac 1\epsilon
$$
thus, it is $-\infty$,  and this fact reflects, when $\epsilon\to 0$,  the density of the zeros of the zeta function.\vspace{.05in}

We emphasize that the role of the Riemann-Weil explicit analytic formulas in the process of overcoming  the initial difficulty through a solution defined by a positive distribution $N(q)$,  directly connects the original (classical geometric) viewpoint with  the trace formula in \cite{Co-zeta}, thus  providing a first  geometric description for  the points of $\bf C$ in terms of the double quotient
\begin{equation}\label{doublequot}
X_\Q:=\Q^\times\backslash \A_\Q/\hatz
\end{equation}
of the adele class space of the rationals by the maximal compact subgroup $\hatz$ of the idele class group. The main key player in this construction is  the scaling action of $\R_+^\times$ which provides\footnote{To remove the divergent logarithmic term from the trace formula \cite{Co-zeta} one needs to remove from $X_\Q$ the orbit of the unit adele $1$, \ie equivalently to subtract the regular representation of $\R_+^\times$ as in \cite{Meyer}.} the above counting distribution $N(u)$, $u\in [1,\infty)$, that determines, in turn, the complete Riemann zeta function via a limiting procedure as $q\to 1$, operated on the Hasse-Weil formula.
Noncommutative geometry plays a crucial role in this development mainly by implementing the noncommutative space $X_\Q$ which naturally arises as the dual of the BC-system \cite{BC}. \vspace{.05in}

To achieve a more classical geometric understanding of the adele class space $X_\Q$ with its scaling action,  in analogy with the action of the Frobenius automorphism on the points of a curve over the algebraic closure of a ground field, one needs to push further the search of other unexpected interactions...
This  geometric understanding  comes in fact from the interplay among three a priori unrelated theories
\begin{enumerate}
\item{Noncommutative Geometry}
\item{Grothendieck topoi}
\item{Tropical Geometry.}
\end{enumerate}
The natural starting point is the  topos  $\wnt$, defined in \cite{CCas} as the Grothendieck presheaf topos dual to the multiplicative monoid $\nt$ of non-zero positive integers. This space is in fact the geometric incarnation of   $\nt$-actions on sets. These actions often occur in global instances of  Frobenius endomorphisms:  for  $\lambda$-rings  they were advocated in \cite{Borgers} in the context of varieties over $\F_1$ ("Futurism" in Manin's interpretation, \cite{M}). Special $\lambda$-rings $R$  (\cite{Atiyah1} Proposition 5.2), belong naturally to the topos $\wnt$ since the Adams operations $\psi_n$ turn $R$ into a sheaf of rings on the topos $\wnt$.  \newline
 At a very basic algebraic level, a fundamental example of Frobenius action of $\nt$ occurs in the theory of semirings (\ie when one drops the existence of the additive inverse in rings). For a semifield\footnote{A semifield is a semiring whose non-zero elements form a group under multiplication} $R$ of ``characteristic one'' (aka idempotent: \ie such that $1+1=1$), the map $x\mapsto x^n=\fr_n(x)$ is an injective endomorphism \cite{Golan}, for any integer $n\in \nt$.
 Thus, one obtains a canonical action of the semigroup $\nt$ on any such $R$. For this reason it is natural to work with the topos $\wnt$  endowed with an action of $\nt$. Furthermore, one also knows that there is a unique semifield\footnote{As a multiplicative monoid $\zmax$ is obtained by adjoining the zero element $-\infty$ to the infinite cyclic group $\Z$ while the operation which plays the role of addition in the semifield is $(x,y)\mapsto \max(x,y) $} $\zmax$ whose multiplicative group is infinite cyclic and it is of characteristic one. Given these facts, it is natural to introduce the following space
 
\begin{definition}[\cite{CCas}]\label{site} The Arithmetic Site $\aarith=\arith$ is the topos $\wnt$
endowed with the  {\it structure sheaf} $\cO:=\zmax$, viewed as a semiring in the topos {\it and} with the action of $\nt$ by  Frobenius endomorphisms.
\end{definition}

\no The semifield $\zmax$ and its companion $\rmax$ (whose multiplicative group is $\R_+^*$), are familiar objects in tropical geometry where the maximum substitutes the usual addition.\newline
By implementing a straightforward generalization in semi-ringed toposes of the understanding of a point  in algebraic geometry, one  obtains the following result which determines a bridge connecting noncommutative geometry with (Grothendieck) topos theory
\begin{theorem}[\cite{CCas}]	
 \label{structure3z} The set of points of the arithmetic site  $\aarith$ over $\rmax$ is canonically isomorphic to $X_\Q=\Q^\times\backslash \A_\Q/\hatz$. The action of the Frobenius
automorphisms $\fr_\lambda$ of $\rmax$ on these points corresponds to the action of the idele class group on $X_\Q=\Q^\times\backslash \A_\Q/\hatz$.
\end{theorem}
\no This theorem sheds new light on a geometric intuition of the curve $\bf C$, in particular, it displays the noncommutative space $X_\Q$  as the set of points of $\bf C$ over the semifield $\rmax$, with the scaling action understood as the action of the Galois group $\Aut_\B(\rmax)$ of $\rmax$ over the Boolean semifield\footnote{$\B:=\{0,1\}$ with $1+1=1$} $\B$. It also suggests that  $\rmax$ ought to be involved in the construction of the  ``algebraic closure'' of $\F_1$, and that the combinatorial core underlying $\bf C$ is countable since both $\nt$ and $\zmax$ are so. We find quite remarkable that while the Arithmetic Site is a combinatorial object of countable nature, it comes nonetheless endowed with a one-parameter semigroup of ``correspondences'' which can be viewed as congruences on the square of this site \cite{CCas}.\vspace{.03in}

The countable set of places of $\Q$ (the points of the Arakelov compactification $\spzb$), is the (classically) visible analog of the set of the orbits of the Frobenius automorphism in the function field case. One obtains a better view of the points of $\bf C$ by considering the periodic orbits $C_p$ (parameterized by primes $p$) as they occur among the points of the Arithmetic Site  $\aarith$ over $\rmax$. One shows that the points of $C_p$ form a circle  whose elements are rank-one subgroups of the multiplicative group of $\rmax$ of the form
\begin{equation}\label{pointdefnp}
H_\mu:=\{\mu^{\frac{n}{p^k}}\mid n\in \Z,\ k\in \N\}.
\end{equation}
This subgroup is unchanged if one replaces $\mu$ with $\mu^p$, and the Frobenius action of $\Aut_\B(\rmax)=\R_+^*$, $\mu\mapsto \mu^\lambda$, induces the transitive action of the quotient group $\R_+^*/p^\Z$. The length of this periodic orbit is $\log p$, and their full collection plays a key role in the trace formula interpretation of the Riemann-Weil explicit formulas in \cite{Co-zeta}. Moreover, each  $C_p$ inherits, as a subspace of the Scaling Site (obtained from the Arithmetic Site by extension of scalars),  a structure sheaf (of characteristic one) which turns each periodic orbit into the analog of a classical elliptic curve \cite{CCscal2}. In this way, one can still apply several key tools of algebraic geometry, such as rational functions, divisors, etc. A striking new feature of the geometry of a periodic orbit is that the degree of a divisor is a real number. For any divisor $D$ in $C_p$, there is a corresponding Riemann-Roch problem with solution space $H^0(D)$.   The continuous dimension\footnote{In analogy with von-Neumann's continuous dimensions of the theory of type II factors} $\cdim(H^0(D))$ of this $\rmax$-module is defined by the limit
\begin{equation}\label{rr1}
\cdim(H^0(D)):=\lim_{n\to \infty} p^{-n}\tdim(H^0(D)^{p^n})
\end{equation}
where $H^0(D)^{p^n}$  is a naturally defined filtration and $\tdim(\cE)$ denotes  the topological dimension of an $\rmax$-module $\cE$.  The following Riemann-Roch formula holds 
\vsp
\begin{theorem}[\cite{CCscal2}]\label{RRperiodic}
$(i)$~Let $D\in \div(C_p)$ be a divisor with $\deg(D)\geq 0$. Then the limit in \eqref{rr1} converges and one has  
\[\cdim(H^0(D))=\deg(D).\]
$(ii)$~The following Riemann-Roch formula holds
\begin{equation*}
\cdim(H^0(D))-\cdim(H^0(-D))=\deg(D)\qquad\forall  D\in \div(C_p).
\end{equation*}
\end{theorem}
In view of these results and the leading role played by the Boolean semifield $\B$ among algebraic idempotent structures\footnote{$\B$ is, in particular, the only finite semifield that is not a field \cf~\cite{Golan}}, one might be (wrongly) induced to think of $\B$ as the natural incarnation of $\F_1$. However, this cannot be the case for the straightforward reason that\footnote{algebras over $\B$ are of characteristic one}: \vspace{.05in}

\centerline{{\it The ring $\Z$ is not an algebra over $\B$.}}

\section{The absolute coefficients, spectra and  $\sss$.}\label{sectcoeffs}

The above undeniable fact led us, once again, to compare  Manin's ideas on  $\F_1$  with another a priori unrelated topic: this is  the world of homotopy theory spectra. Topological spectra greatly generalize cohomology theories; many important invariants in algebraic topology, like ordinary cohomology and K-theory, can be reformulated in terms of spectra, which thus provide a unified treatment for ``generalized coefficients''. One fundamental discovery in the topological context is that ``ring spectra'' generalize rings, and in particular, the ``sphere spectrum'' $\underline\sss$ becomes more basic than the ring $\Z$, because the latter can be seen as an algebra over the former. This theory of ``brave new rings'' has proved to be the right framework for cyclic homology; in particular,  the theory of $\Gamma$-spaces is known to provide a workable model of connective spectra \cite{DGM}. One usually works at the homotopy level, so it is crucial to handle Kan complexes to obtain a good model structure. However, to take full advantage of this theory for the development of  Manin's ideas on $\F_1$ in number theory, we believe that  $\Gamma$-spaces ought to be viewed in their most basic form, namely  as simplicial objects in the category of $\Gamma$-sets, so that homotopy theory can play the role of homological algebra corresponding to the ``absolute algebra''  over the base $\Gamma$-ring $\sss$  \cite{CCprel}. This $\Gamma$-ring  is the categorical starting point in the construction of the  sphere spectrum $\underline\sss$, together with the natural functor from $\Gamma$-spaces to spectra, and its is exactly this basic combinatorial nature that makes it closer to the sought for $\F_1$. The category $\gam$ of pointed $\Gamma$-sets (aka $\sss$-modules $\catmo(\sss)$) can be directly described as follows.  One starts with the small category $\Gamma^{\rm op}$  as a full subcategory of the category of finite pointed sets whose objects are the pointed finite sets\footnote{where $0$ is the base point.} $k_+:=\{0,\ldots ,k\}$, for $k\geq 0$.  In particular, $0_+$ is both initial and final in $\Gamma^{\rm op}$,   making $\gop$ a  {\it pointed category}. A $\Gamma$-set is defined as a (covariant) functor $\gop\longrightarrow\Ses$ between pointed categories, and 
the morphisms in this category are natural transformations. One lets $
\sss:\gop\longrightarrow \Ses
$
be the inclusion functor. The internal hom functor is defined by
$$\underline\Hom_\sss(M,N):=\{k_+\mapsto \Hom_\sss(M,N(k_+\wedge -))\}.
$$
This formula uniquely defines the smash product of $\Gamma$-sets by applying the adjunction 
$$
\underline\Hom_\sss(M_1\wedge M_2,N)=\underline\Hom_\sss(M_1,\underline\Hom_\sss(M_2,N)).
$$
The basic construction of $\sss$-modules  associates to an abelian monoid $A$ with a zero element, the Eilenberg-MacLane functor $M=HA$  
$$
HA(k_+)=A^k, \qquad  Hf:HA(k_+)\to HA(n_+), $$ $$ \ Hf(m)(j):=\sum_{f(\ell)=j} m_\ell,
$$
where $m=(m_1,\ldots ,m_k)\in HA(k_+)$, and the zero element of $A$ gives meaning to the empty sum. An $\sss$-algebra $\mathcal A$  is an $\sss$-module $\mathcal A: \gop\longrightarrow\Ses$  endowed with an associative multiplication 
$\mu:\cA \wedge \cA\to \cA$ and a unit~ $1: \sss\to \cA$.\newline An ordinary semiring $R$ gives rise to the $\sss$-algebra $HR$, and the corresponding embedding of categories is fully faithful so that no information is lost. In contrast, the basic $\sss$-algebra $\sss$  now lies under $HR$ for any semiring $R$.\newline
Given a multiplicative monoid $M$ with a zero element $0\in M$ such that $0\times x=x\times 0=0$ for all $x\in M$, one defines the spherical $\sss$-algebra $\sss[M]$ which associates to the pointed set $X$ the smash product $X\wedge M$, where the base point of $M$ is $0\in M$. One identifies $\sss[M][1_+]=1_+\wedge M$
with $M$ by sending the base point of $1_+\wedge M$ to $0\in M$, and $a\wedge m$ where $a\in 1_+\setminus \{*\}$ and $m\in M\setminus \{0\}$ to $m$. To avoid confusion we write $2_+=\{*,a,b\}$. Besides the base point the elements of $\sss[M][2_+]=2_+\wedge M$ are given by pairs of the form $(a,m)$ or $(b,m)$ where $m\in M\setminus \{0\}$. 
One has three natural pointed maps $f:2_+\to 1_+$, which are
$$
\phi(a)=a,\ \phi(b)=*, \ \ \psi(a)=*, \ \psi(b)=a, \ \ \rho(a)=\rho(b)=a.
$$
Let $m\in M\setminus \{0\}$ and consider the pair $z=(b,m)\in \sss[M][2_+]$. One has
$\phi_*(z)=*=0$ and $\psi_*(z)=m$. Moreover one has $\rho_*(z)=m$. This means that for the partially defined addition in $\sss[M][1_+]=M$, one has $0+m=m$ for all $m\in M$.\newline
Thus both ordinary rings and monoids fit fully faithfully and naturally \cite{CCprel} (Proposition 3.5), in the category of $\sss$-algebras  yielding a strong argument for viewing $\sss$ as the natural candidate for $\F_1$. Nonetheless one needs to test this idea in various ways. For instance, one sees  \opcit that the tensor square of $H\Z$ over $\sss$ is non-isomorphic to $H\Z$,  and this result provides more ground to the original intuition of Manin in \cite{Man-zetas}. One may also wonder which advancements this point of view may produce to the understanding of the ring $\Z$ and its algebraic spectrum $\Spec \Z$. We shall now move to a detailed discussion of this topic.  \vspace{.03in}

Let $\spzb$ be the Arakelov compactification of  $\Spec \Z$ obtained by adding the archimedean place with associated symbol  $\infty$. Then, the new point of view described above provides a natural extension of the classical structure sheaf of $\Spec \Z$ to the Arakelov compactification. The crucial points concerning the quest for the curve $\bf C$ are two: firstly,  this extended structure sheaf $\cO$ is still a subsheaf of the constant sheaf $\Q$; the second interesting point is that the global sections of $\cO$ form a finite algebra extension of $\sss$. This extension is identifiable with the extension by the two roots of unity inside $\Q$ that we used in  \cite{CCn1} in the process of showing that Chevalley groups are varieties over $\F_{1^2}$ in the sense of Soul\' e\footnote{Another convincing argument in favor of $\sss$-algebras is that the ad-hoc category we introduced in \cite{CC0} to simplify Soul\' e's definition of varieties over $\F_1$, is naturally (see \cite{schemeF1}) a full subcategory of the category of $\sss$-algebras}. The condition that restricts the elements of $H\Q$ at the archimedean place is simple to formulate when one views the functor $H\Q$ as assigning to a finite pointed set $F$ the $\Q$-valued divisors on $F$. The restriction is then stated by writing that the sum of the absolute values of the involved rational numbers is $\leq 1$. One checks that this condition is stable under push-forwards and products and hence it defines a sub-$\sss$-algebra of $H\Q$. This sub-$\sss$-algebras, defined using a norm, also applies in the context of the adeles of a global field and allows one to transpose the approach due to A. Weil of the Riemann-Roch theorem for function fields to the number field $\Q$ \cite{RR}. 

\no A divisor $D$ on $\spzb$ defines a compact {\it subset} $K=\prod K_v\subset \A_\Q$ of the locally compact ring of adeles. When $p$ is a non-archimedean prime, each $K_p\subset \Q_p$ is an additive subgroup, in contrast, the compact subset $K_\infty \subset \R$ is just a symmetric interval whose lack of additive structure prevents one to use Weil's original construction involving the addition map $\psi:\Q \times K \to \A_\Q$. On the other hand, one also quickly notices that  $\psi$ retains its meaning in the context of $\sss$-modules,  giving rise to a short complex. Using the Dold-Kan correspondence in the context of $\sss$-algebras, one then introduces a $\Gamma$-space $H(D)$ which encodes the homological information of the divisor $D$ and only depends upon the linear equivalence class of $D$ (\ie the divisor class is unchanged under the multiplicative action of $\Q^\times$ on $\A_\Q$). As a by-product, one obtains a Riemann-Roch formula for Arakelov divisors on $\spzb$ of an entirely novel nature that relies on the introduction of three new key notions: (integer) dimension, cohomologies $(H^0(D), H^1(D))$ (attached to a divisor $D$), and Serre duality. More precisely, the Riemann-Roch formula equates the integer-valued Euler characteristic with a simple modification of the traditional expression (\ie the degree of the divisor plus log 2).   
\begin{theorem}[\cite{RR}] \label{rrspzbintro} 
Let $D$ be an Arakelov divisor  on $\spzb$. Then 
\begin{equation}\label{rrforq1intro}
\dim_{\spm}H^0(D)-\dim_{\spm}H^1(D)=\bigg\lceil \deg_3 D+ \log_3 2\bigg\rceil	-\mathbf 1_L.
\end{equation}
Here,  $\lceil x \rceil$ denotes the odd function on $\R$ that agrees with the ceiling function on positive reals, and $\mathbf{1}_L$ is the characteristic function of an exceptional set\footnote{$L\subset \R$ is the union, for $k\geq 0$, of the intervals $\deg(D)\in(\log \frac{3^k}{2},\log \frac{3^k+1}{2})$} of finite Lebesgue measure.
\end{theorem}
\no In  \eqref{rrforq1intro},  the neperian logarithm that is traditionally used to define the degree of a divisor $D=\sum_j a_j\{p_j\} + a\{\infty\}$ in Arakelov geometry, is replaced by the logarithm in base $3$. This alteration is equivalent to the division by $\log 3$ i.e.  $\deg_3(D):=\deg(D)/\log 3$, $\log_3 2=\log 2/\log 3$.\newline
The number $3$ appears unexpectedly in the computation of the dimension of the cohomology of the $\spm$-modules  by determining their minimal number of linear generators. For 
 $\dim_{\spm}H^0(D)$ one finds that the most economical way of writing the elements of a symmetric interval $\Z\cap K_\infty$ involves writing integers as polynomials of the form 
 \begin{equation}\label{zaspol}
 \sum_{j\ge 0} \alpha_j \ 3^j, \ \ \alpha_j\in \{-1,0,1\}.
\end{equation}
 Similarly, in the case of $\dim_{\spm}H^1(D)$,  one finds that the best way to approximate elements of the circle $\R/\Z$ is to use Laurent polynomials of the form 
 \begin{equation}\label{zaspol1}
 \sum_{j< 0} \alpha_j \ 3^j, \ \ \alpha_j\in \{-1,0,1\}.
 \end{equation}
 The key fact here is that the addition\footnote{once the addition is defined, the product follows uniquely using $X^j\, X^k=X^{j+k}$} of polynomials $P(X)=\sum_{j\geq 0} \alpha_j \ X^j, \ \ \alpha_j\in \{-1,0,1\}$ with coefficients in $\spm$ is 
identical to the addition of (truncated) Witt vectors for the finite field $\F_3$. One finds that the addition $P+Q$ of two polynomials of degree $\leq n$ gives a polynomial of degree $\leq n+1$, and that the only non-obvious rule one has to prescribe is the sum: $1+1:= X-1$. Conceptually, the fundamental point is that the image of the Teichmuller lift for $\F_3$ sits inside $\Z$. At the same time, the Witt vectors with only finitely many non-zero components form a subring of the Witt ring, and this subring is  $\Z$!   

\section{The ring of integers as a ring of polynomials}\label{sec4}

There is another way to represent the integers as polynomials in one variable, and in this description, the ``coefficients'' belong to the absolute base $\sss$. This representation is known as the {\em negabinary} representation of numbers  
\begin{equation}\label{negabinary}
n=\sum \alpha_j\ (-2)^j, \ \ \alpha_j\in \{0,1\}.
\end{equation} 
The number $X=-2$ is remarkably unique, making the representation of an integer $n$ possible as polynomial $P(X)$ with coefficients $\alpha_j\in \{0,1\}$. By following the same steps that led us to  Theorem \ref{rrspzbintro}, but working now over the absolute base $\sss$, one obtains the following new and simplified version of the Riemann-Roch formula which now involves the logarithm in base $2$ 
\begin{theorem}[\cite{RR1}]\label{rrz} 
Let $D$ be an Arakelov divisor  on $\spzb$. Then 
\begin{equation}\label{rrforq1}
\dim_{\sss}H^0(D)-\dim_{\sss}H^1(D)=\bigg\lceil \deg_2 D\bigg\rceil'	+1
\end{equation}
where $\lceil x \rceil'$ is the right continuous function which agrees with ceiling$(x)$ for $x>0$ non-integer and with $-$ceiling$(-x)$ for $x<0$ non-integer.
\end{theorem}
\no This  version of the Riemann-Roch Theorem improves on Theorem \ref{rrspzbintro} for the following reasons:
\begin{enumerate}
\item The term $\mathbf 1_L$ involving the exceptional set $L$ in the original statement (see \cite{RR}) has now disappeared from the formula.
\item The formula \eqref{rrforq1} is in perfect analogy with the Riemann-Roch theorem for curves of genus $0$.
\item The canonical divisor  $K=-2\{2\}$ has integral degree $\deg_2(K)=-2$.
\end{enumerate}
Theorem~\ref{rrz}  fits now perfectly with the tri-lingual text suggested by A. Weil, that supports the analogy between Riemann's transcendental theory of algebraic functions of one variable in the first column, the algebraic geometry of curves over finite fields, in the middle column, and the theory of algebraic number fields in the third column. Indeed, according to Weil
\begin{quote}
\emph{Mais on peut, je crois, en donner une idée imagée en disant que le mathématicien qui étudie ces problèmes, a l'impression de déchiffrer une inscription trilingue. Dans la première colonne se trouve la théorie riemannienne des fonctions algébriques au sens classique. La troisième colonne c'est la théorie arithmétique des nombres algébriques. La colonne du milieu est celle dont la découverte est la plus récente : elle contient la théorie des fonctions algébriques sur un corps de Galois. Ces textes sont l'unique source de nos connaissances sur les langues dans lesquels  ils sont \' ecrits; de chaque colonne, nous n'avons bien entendu que des fragments ; la plus compl\`ete et celle que nous lisons le mieux, encore \`a pr\' esent, c'est la premi\`ere. Nous savons qu'il y a de grandes diff\' erences de sens d'une colonne \`a l'autre, mais rien ne nous en avertit \`a l'avance. \'A l'usage, on se fait des bouts de dictionnaire, qui permettent de passer assez souvent d'une colonne à la colonne voisine.}
\end{quote}

\no In  Weil's vision there is, in the middle column (that of function fields),  a geometric understanding of the zeta function as the generating function of the number of points of the curve over extensions of the field of constants. In section~\ref{cincarnation} we translated in this dictionary the Hasse-Weil formula, thus leading one to the first encounter with the ``the curve'' $\bf C$ and the action of $\R^*_+$ on $\bf C$, analogous to a Galois action. Theorem~\ref{rrz} indicates that the role of the field of constants is played by the absolute coefficient ring $\sss$. Since the boolean semifield $\B$ can be viewed as a $\sss$-algebra, this translation suggests to descend the structures of the Arithmetic and Scaling Sites discussed in section~\ref{cincarnation} from $\B$ to $\sss$. 

\section{Rings of \texorpdfstring{$\sss[\mu_{n,+}]$}--polynomials}\label{sec5}

Let $n>0$ be an integer, $\mu_n$ the multiplicative group of $n$-th roots of $1$ and  $\sss[\mu_{n,+}]$ the spherical $\sss$-algebra of the (pointed) monoid $\mu_{n,+}=\mu_n\cup \{0\}$. We recall that morphisms of $\sss$-algebras $\sss[\mu_{n,+}]\to HR$ correspond (bijectively) to group homomorphisms $\iota: \mu_n\to R^\times$ \cite{CCprel}. In this section, we introduce the notion of rings of $\sss[\mu_{n,+}]$-polynomials in one (Definition \ref{sgenerator})  and several variables (Remark~\ref{frame})  which might play a key role in the search of the ``absolute Descartes powers''  among ordinary rings. We show that the pair $(R,X)$ of a ring $R$ and an $\sss[\mu_{n,+}]$-generator of $R$ is uniquely characterized, up to isomorphism, by the map from $\mu_n$ to polynomials with coefficients in the pointed monoid $\mu_{n,+}$, that encodes the addition of $1$ into elements of $\mu_n$.

\begin{definition} \label{sgenerator} Let $R$ be a ring,  $\iota:\mu_n\to R^\times$ be an injective group homomorphism. An element $X\in R$ is an $\sss[\mu_{n,+}]$-generator of $R$ if and only if every element  $z\in R$ can be  written uniquely as a polynomial $z=\sum_j \iota(\alpha_j) \, X^j$ with coefficients $\alpha_j\in \mu_n\cup \{0\}$.
\end{definition}
\begin{remark}\label{frame} More generally,  a finite set $\{X_i\mid i\in \{1,\ldots ,k\}\}$, $\sss[\mu_{n,+}]$-generates $R$ if and only if every element  $z\in R$ can be written  uniquely as a polynomial $z=\sum_j \iota(\alpha_j) \, X^j$ with coefficients $\alpha_j\in \mu_n\cup \{0\}$, where $j$ is  a multi-index $j=(j_1,\ldots, j_k)\in \N^k$, and $X^j:=\prod X_i^{j_i}$.    
\end{remark} 

Let $\cP(\mu_n)$ be the subset of the set $\left(\mu_n\cup \{0\}\right)^\N$ of sequences with only finitely many non-zero terms.  Let $X\in R$, then the map $
\sigma:\cP(\mu_n)\to R$, given by
\begin{equation}\label{sigma}
\sigma((\alpha_j)):=\sum_j \iota(\alpha_j) \, X^j
\end{equation}   is well defined since for $\alpha=(\alpha_j)\in \cP(\mu_n)$ the sum $\sum_j \iota(\alpha_j) \, X^j$ defines an element of $R$. It follows from  Definition \ref{sgenerator} that if $X$ is an $\sss[\mu_{n,+}]$-generator, the map $\sigma$ is a bijection of $\cP(\mu_n)$ with $R$.
\vspace{.05in}

\no The simplest instance of a $\sss[\mu_{n,+}]$ generator, with $n+1$ a prime power $q$, is provided by the following example. 

\begin{example}\label{sexample1} The ring $R=\F_q[X]$ of polynomials over the finite field  $\F_q$ has the variable $X$ as  $\F_q^\times$-generator. 	
\end{example}

\no Next proposition shows that the $m$-th root  of an $\sss[\mu_{n,+}]$-generator $X$ of a ring $R$ is a $\sss[\mu_{n,+}]$-generator of the $R$-algebra extension $R[Y]/(Y^m-X)$, hence providing an infinite source of examples. 
\begin{proposition}\label{sprop1} 
Let $R$ be a ring,  $\iota:\mu_n\to R^\times$ be an injective group homomorphism, $X\in R$  an $\sss[\mu_{n,+}]$-generator of $R$ and $m\in \N$ be a positive integer. Then $Y\in R[Y]/(Y^m-X)$ is an $\sss[\mu_{n,+}]$-generator of $R[Y]/(Y^m-X)$.	
\end{proposition}
\proof 
Any element  $z$ of $R[Y]/(Y^m-X)$ can be  written uniquely as $z=\sum_{j=0}^{m-1} a_j Y^j$, with $a_j\in R$  written uniquely as $a_j=\sum_{j,k} \iota(\alpha_{j,k}) \, X^k$
where $\alpha_{j,k}\in \mu_n\cup \{0\}$. Since $Y^m=X$ one obtains the unique finite decomposition 
$$
z=\sum_{j,k} \iota(\alpha_{j,k}) Y^{j+m k}, \qquad \alpha_{j,k}\in \mu_n\cup \{0\}.
$$
\endproof 
The following example is a straightforward generalization of the fact that $3$ is an $\spm=\sss[\mu_{2,+}]$-generator of the ring $\Z$ of integers.

\begin{example}\label{sxam1}
	Let $m\in \N$ be a positive integer, and $\epsilon =\pm 1$. Then $X=(3\epsilon)^{1/m}$ is an $\spm$-generator of the subring $R=\Z[X]$ of the number field $\Q((3\epsilon)^{1/m})$. 
 
 \no Indeed, the polynomial $X^m-3\epsilon$ is irreducible, thus every element  $z\in R$ can be  written uniquely as a sum
$$
z=\sum_{j=0}^{m-1} a_j X^j, \qquad a_j\in \Z.
$$
In turns, every $a_j$ can be uniquely written as $a_j=\sum_{j,k} \alpha_{j,k} \, (3\epsilon)^k$,
where $\alpha_{j,k}\in \{-1,0,1\}$. Since $3\epsilon=X^m$ one obtains the unique decomposition 
$$
z=\sum_{j,k} \alpha_{j,k} X^{j+m k}, \qquad \alpha_{j,k}\in \{-1,0,1\}.
$$
An interesting case is for $m=2$ and $\epsilon=-1$ since then the ring $R=\Z[\sqrt{-3}]$  is an order of the ring of integers of the imaginary quadratic field $\Q(\sqrt{-3})$. \end{example}
Notice that in the Example~\ref{sxam1} the addition is specified by an equality of the following form 
\begin{equation}\label{polX}
1+1=P(X), \qquad 	P(X)=\sum_j \alpha_j \, X^j, \qquad \alpha_j\in \{-1,0,1\},
\end{equation}
with   $P(X)=\epsilon\,X^m-1$. A simple algebraic presentation of the form \eqref{polX} holds when working over $\mu_{n,+}$ for $n=1,2$.\vspace{.05in}

\no  The following result states the uniqueness of  a similar polynomial presentation in the general case.

\begin{proposition}\label{sprop2} Let $R$ be a ring, $\iota:\mu_n\to R^\times$ be an injective group homomorphism,  $X\in R$  an $\sss[\mu_{n,+}]$-generator of $R$. For a polynomial decomposition $z=\sum_j \iota(\alpha_j) \, X^j\in R$, let $\deg(z)$ be the smallest integer $m$ such that $\alpha_j=0$ for all $j>m$. Then, the following results hold
\begin{enumerate}
\item[(i)] Let $m\in \N$, and $J_m=\langle X^m\rangle \subset R$ be the ideal generated by $X^m$. Any element $z\in R$ admits a unique decomposition as $z=a+b$ where $\deg(a)<m$ and $b\in J_m$.
\item[(ii)] The quotient $R_m:=R/J_m$ is a finite ring whose elements are uniquely written as $\sum_{j=0}^{m-1} \iota(\alpha_j) \, X^j$, with $\alpha_j\in \mu_{n,+}=\mu_n\cup \{0\}$.
\item[(iii)] The quotient $R_1:=R/J_1$ is a finite field with $n+1$ elements and $\iota:\mu_{n,+}\to R$ is a multiplicative section of the quotient map $R\to R_1$. 
\item[(iv)] The canonical ring homomorphism $\pi:R\to \varprojlim R_m$ is injective.
\item[(v)] The pair $(R,X)$ is uniquely specified, up to isomorphism, by the map $h:\mu_n\to \cP(\mu_n)$ which is uniquely defined by the equality $\sigma(h(\xi))=\iota(\xi)+1$.
\end{enumerate}
\end{proposition}

\proof $(i)$~Let $z=\sum_j \iota(\alpha_j) \, X^j$. By writing $z$ as 
\begin{equation}\label{decab}
z=\sum_{j=0}^{m-1} \iota(\alpha_j) \, X^j+\sum_{j=m}^{\deg(z)} \iota(\alpha_j) \, X^j=a +X^m c
\end{equation}
one obtains the required decomposition with $b=X^m\, c$. The uniqueness of such decomposition then follows 
from the uniqueness of the decomposition as in Definition \ref{sgenerator}.\newline 
$(ii)$~Follows from $(i)$.  In particular, one easily checks that $R_m$ has cardinality $\#(R_m)=(n+1)^m$.
\newline
$(iii)$~By construction the map $\iota:\mu_{n,+}\to R$ is a multiplicative section of the quotient map $R\to R_1$. It follows that the non-zero elements of $R_1$ form the multiplicative group $\mu_n$ so that $R_1$ is a field with $n+1$ elements.\newline
$(iv)$~The components of $z=\sum_j \iota(\alpha_j) \, X^j\in R$ are uniquely determined by $\pi(x)$.
\newline
$(v)$~Let $(R',X')$ be a second pair corresponding to the same map $h:\mu_n\to \cP(\mu_n)$. Let $\rho:R\to R'$ be the bijective map defined by 
$$
\rho\Big(\sum_j \iota(\alpha_j)\, X^j\Big):=\sum_j \iota'(\alpha_j) \, X'^j, \qquad \alpha_j\in \mu_n\cup \{0\}.
$$
 One has by construction
\begin{equation}\label{polX1}
\deg(a)<m~ \Longrightarrow~ \rho(a+X^mb)=\rho(a)+(X')^m\rho(b), \qquad \forall b.
\end{equation}
 In particular one also has $\rho(J_m)=J'_m$ for all $m$. Thus $\rho$ induces a bijection $\rho_m:R_m\to R'_m$. By $(iii)$, to show that $\rho$ is a ring homomorphism, it is enough to verify that each $\rho_m$ is a ring homomorphism. \newline
 To show that $\rho_m$ is additive it is enough to show that one can compute all the components of a sum 
 \begin{equation}\label{summ}
 \sum_{j=0}^{m-1} \alpha_j \, X^j+\sum_{j=0}^{m-1} \beta_j \, X^j=
\sum_{j=0}^{m-1} \gamma_j \, X^j
 \end{equation}
using only the map $h:\mu_n\to \cP(\mu_n)$. To do this one first determines a map $F$ from $k$-tuples of elements of elements of $\mu_{n,+}$ to pairs $(x,Z)$ where $x\in \mu_{n,+}$ and where $Z$ is a $(k-1)$-tuple of elements of $\cP(\mu_n)$.
The map $h$ determines uniquely a  symmetric map 
\[
H:\mu_{n,+}\times \mu_{n,+}\to \mu_{n,+}\times \cP(\mu_n), \ \ H(\xi,\eta)=(\xi+\eta,0) \ \ \text{if}\ \ \xi\ \eta=0
\]
\begin{equation}\label{defnH}
H(\xi,\eta)=(H_0(\xi,\eta),P(\xi,\eta)), \ \ H_0(\xi,\eta)+X P(\xi,\eta)=\eta\ h(\xi\  \eta^{-1})\ \text{if}\ \ \eta \neq 0
\end{equation}
 To define $F$ one proceeds by induction on $k$. For $k=1$ one lets $F(x)=x$. For $k=2$ one lets $F_2=H$.
We denote the two components of  
$
F_k:\mu_{n,+}^{k-1}\times \mu_{n,+}\to \mu_{n,+}\times \cP(\mu_n)^{k-1}
$ as $F_k^{(1)}$ and $F_k^{(2)}$. To pass from $k-1$ to $k$ one lets
$$
F^{(1)}_k(\alpha,\eta):=(H_0(F^{(1)}_{k-1}(\alpha),\eta),\ \ F_k^{(2)}(\alpha,\eta):= (F^{(2)}_{k-1}(\alpha),P(F^{(1)}_{k-1}(\alpha),\eta))
$$
where in the last expression we append the polynomial $P(F^{(1)}_{k-1}(\alpha),\eta)$ to the list $F^{(2)}_{k-1}(\alpha)$, thus obtaining a list of $k-1$ polynomials.\newline
 To compute the components $\gamma_j$ of the sum \eqref{summ}, we build by induction on $k$, two lists. The first $R(k)$ is the list of the coefficients already computed and it is the single list given by $(\gamma_0,\gamma_1,\ldots,\gamma_{k-1})$. The second $C(k)$, (called the carry), is a list of polynomials with coefficients in $\mu_{n,+}$ and it is encoded as the list of their coefficients. Each such list $\ell$ of coefficients has $m-k$ terms, all in $\mu_{n,+}$. We denote by $f(\ell)\in \mu_{n,+}$ the first term of the list $\ell$ and by $t(\ell)$ the list obtained by dropping the first element of the list $\ell$, it has $m-k-1$ terms.
The step to obtain $R(k+1),C(k+1)$ from $\alpha,\beta,R(k),C(k)$ is 
$$
R(k+1):=F^{(1)}(\alpha_k,\beta_k,(f(\ell))_{\ell\in C(k)})
$$
and 
$$
C(k+1):=(t(\ell))_{\ell\in C(k)}, F^{(2)}(\alpha_k,\beta_k,(f(\ell))_{\ell\in C(k)})
$$
where one replaces each element of  $F^{(2)}(\alpha_k,\beta_k,(f(\ell))_{\ell\in C(k)})$ by the list of  its first $m-k$ coefficients.\newline 
More concretely 
    one first obtains 
$
\gamma_0=F^{(1)}_2(\alpha_0,\beta_0)
$
while the carry over delivers the polynomial 
$
P(\alpha_0,\beta_0)=F^{(2)}_2(\alpha_0,\beta_0)
$. Thus $R(1)=(\gamma_0)$, $C(1)$ has one element which is the list of the first $m-1$ coefficients of $P(\alpha_0,\beta_0)$. One then trims the elements $\alpha, \beta$ to and considers the sum 
\begin{equation}\label{summ1}
 \sum_{j=1}^{m-1} \alpha_j \, X^j+\sum_{j=1}^{m-1} \beta_j \, X^j+ X P(\alpha_0,\beta_0)
 \end{equation}
 All terms in \eqref{summ1} are divisible by $X$ and one can use $F_3$ to compute the sum of the three terms $\alpha_1,\beta_1,p_0$ where $p_0$ is the constant term of $P(\alpha_0,\beta_0)$. This operation delivers the next term $$\gamma_1=F^{(1)}_3(\alpha_1,\beta_1,p_0)$$ of \eqref{summ}, and adjoins the two polynomials of the list 
$F^{(2)}_3(\alpha_1,\beta_1,p_0)$ to the list of carry over consisting of the single polynomial $P(\alpha_0,\beta_0)$ with its first term $p_0$ deleted. The carry over list consists now of three terms $\ell_1,\ell_2,\ell_3$. One then uses $F_5$ to compute the sum of the $5$ terms : $\alpha_2,\beta_2$ and the three terms $f(\ell_1),f(\ell_2),f(\ell_3)$ from the carry over. This adds $4$ terms to the list of carry over which has now $7$ terms, where the three previous ones have been trimmed by deleting their lowest term. After $k$ such steps the carry over list has $2^k-1$ elements and one proceeds as follows. One uses $F_{2^k+1}$ to compute the sum of the $2^k+1$ terms given by  $\alpha_k,\beta_k$ together with the  terms $f(\ell)$ of the carry over list. This operation delivers $\gamma_k$ and adjoins $2^k$ terms to the carry over list which now consists of $2^{k+1}-1$ terms. This process terminates when $k=m$ and $R(m)$ delivers universal formulas for the terms $\gamma_j$, $0\leq j\leq m-1$ using only $\alpha,\beta$ and the map $h$.\newline 
 The  fact that the coefficients $\gamma_j$ can be computed using only $\alpha,\beta$ and the map $h$ proves that $\rho$ is additive since one can use the same formula to compute $\alpha+\beta$ in $R_m$ and $\rho_m(\alpha)+\rho_m(\beta)$ in $R'_m$. The multiplicativity of $\rho$ follows by bilinearity from $\rho(\alpha X^n\times \beta X^m)=\rho(\alpha X^n)\rho(\beta X^m)$. This shows that $\rho:R\to R'$ is a ring isomorphism  and by construction one has $\rho(X)=X'$.\endproof 
\begin{definition} \label{defnh1}
The map 
\begin{equation}\label{defnh}
   h:\mu_n\to \cP(\mu_n), \qquad  \sigma(h(\xi))=\iota(\xi)+1
\end{equation}
 which characterizes the pair $(R,X)$ (by Proposition~\ref{sprop2}) is called the {\em hold} of the pair $(R,X)$.
\end{definition}

\begin{corollary} Let $n$ be such that there exists a polynomial ring in one generator over $\sss[\mu_{n,+}]$, then $n+1$ is a prime power.
\end{corollary}
\proof This follows from Proposition \ref{sprop2} $(iii)$. \endproof 
 \begin{remark}\label{rem1}
     The proof of 
 Proposition~\ref{sprop2} $(v)$ is stated so that one can, by following it, write a computer program which can be used  to test the additive structure of the ring $R_m$. This will be relevant in  section \ref{sec6} to determine in the various examples the rings $R_m$. \newline
 The map $h:\mu_n\to \cP(\mu_n)$ of \eqref{defnh} determines the addition   $H:\mu_{n,+}\times \mu_{n,+}\to \mu_{n,+}\times \cP(\mu_n)$, \eqref{defnH}, of pairs of elements of $\mu_{n,+}$ using the compatibility with multiplication by elements of $\mu_n$. \newline 
 Proposition \ref{sprop2} shows that a pair $(R,X)$, where $X$ is an $\spm$-generator of $R$, \ie $n=2$, is uniquely characterized by the polynomial $P(X)$ as in  \eqref{polX}. The polynomial $P(X)=-1$ produces the pair $(\F_3[X],X)$,  while $P(X)=X-1$ determines the pair $(\Z,3)$.\newline  When $n=2$, the constant term of the polynomial $P(X)$  in 
\eqref{polX} is necessarily equal to $-1$. Indeed, had the constant term be $0$ or $1$, one would contradict the uniqueness of the  decomposition of Definition \ref{sgenerator} by the equality $1=P(X)-1$. This also shows that $R_1=\F_3$.
\end{remark}

\begin{remark}\label{rem3} It is not true that a random choice of a polynomial with coefficients in $\spm$ and constant coefficient $-1$ corresponds to a pair. A simple case is with $P(X)=-1+X+X^2$. Indeed, in the following lines we show that  $5$ is not represented by any polynomial. With this rule, one has
$
1+1+1+1=1+X+X^2
$. Adding $1$ to both sides gives 
\begin{align*}
&1+1+1+1+1=-1+X+X^2+X+X^2=-1+X(-1+X+X^2)+X^2(-1+X+\\&+X^2)=-1-X+X^3+X^3+X^4=-1-X+X^3(-1+X+X^2)+X^4=\\
&= -1-X-X^3+X^4+X^4+X^5.
\end{align*}
Then, when working in $R_n$ (\ie modulo $X^n$) the number $5$ is represented by 
$$
5=-1-X-X^3-X^4-X^5-\cdots -X^{n-1}\in R_n
$$
and this expression is of degree $n-1$ for any $n$ and thus does not correspond to a finite sum of powers of $X$.
\end{remark}

\vspace{.05in}

\section{Examples}\label{sec6}
In this section we give examples of  polynomial rings $(R,X)$ in one generator $X$ over $\sss[\mu_{n,+}]$ where $R$ is of characteristic zero.  The ring $R$ is embedded as a subring of $\C$ by solving for $X\in \C$ the equations $\sigma(h(\xi))=\iota(\xi)+1, \, \xi \in \mu_n$, using  the canonical embedding  $\mu_{n,+}\subset \C$. The projective limit $\varprojlim R_n$ is, in these examples, a finite extension of the ring of $p$-adic integers $\Z_p$.  While one can use the axiom of choice to show the existence of an embedding of the $p$-adic field $\Q_p$ in the field of complex numbers, such embeddings have the status of a chimera. Indeed, the  continuity of measurable characters of compact groups applied to the additive group $\Z_p$ shows that an embedding of the $p$-adic field $\Q_p$ in the field of complex numbers is automatically non-measurable. On the other hand, next examples will show that polynomial rings $(R,X)$ in one generator $X$ over $\sss[\mu_{n,+}]$ provide instances of explicit interactions of $p$-adic fields (and their finite extensions) with the complex numbers. These interactions are given by pairs of embeddings with dense ranges 
\[
\renewcommand{\arraystretch}{1.3}
\begin{array}[c]{ccccccccc}
\F_q&\stackrel{\pi}{\twoheadleftarrow}&W(\F_q)&\hookrightarrow&W(\F_q)[\eta]&\hookleftarrow& R[X^{-1}]&\hookrightarrow& \C
\end{array}
 \]
of  the ring of Laurent polynomials $R[X^{-1}]$. The left embedding in the above diagram is in a finite algebraic extension $W(\F_q)[\eta]$ of the Witt ring $W(\F_q)$. The field of fractions of the ring $W(\F_q)[\eta]$ is a finite extension of the $p$-adic field. Most of these examples  come from known number systems and have their origin in the search of optimal manners of encoding numbers \cite{Knuth}. In each case, the quotient $R_1=R/(XR)$ is  the finite field $\F_q$, $q=n+1$, and the multiplicative semi-group isomorphism $j:\F_q\sim \mu_{n,+}\subset \C$ serves as a guide, using the addition in the finite field $\F_q$, for the terms of degree $0$ in the construction of the map $h$.  Note that the choice of $j$ for $\bar \F_q$  plays a key role in the construction by Quillen \cite{Q} of the relation between the algebraic $K$-theory of $\F_q$ and the Adams operations.
\subsection{ Polynomial rings  in one generator over $\sss=\sss[\mu_{1,+}]$}
 When working over $\sss=\sss[\mu_{1,+}]$ there is no cancellation since there is no minus sign available. Thus starting from two non-zero elements $x,y$ the equality $x+y=0$ can only be verified in the projective limit $\varprojlim R_m$. We compute this projective limit in the  next examples.
\subsubsection{The polynomial ring $(\Z,-2)$}\label{Z2}
The ring $\Z$ admits the generator $X=-2$ over $\sss$. The hold  is given by $1+1=P(X)=X+X^2$. The values of the polynomials of degree $n$,  at $X=-2$ are reported for the first values of $n$ in the following table
\begin{displaymath}
\begin{array}{|c |c|}
\hline
n & \{p(-2): \deg p= n\}\\ 
\hline  
0  & [0,1]\cap \Z\\
1  & [-2,-1]\cap \Z\\
2 &  [2,5]\cap \Z\\
3 &  [-10,-3]\cap \Z\\
4 & [6,21]\cap \Z\\
5& [-42,-11]\cap \Z\\
6& [22,85]\cap \Z\\
\hline
\end{array}
\end{displaymath}
Let us look, for example, at the computation of $1+1+X$. One gets
$$
1+1+X=X+X^2+X=X(1+1+X)
$$
and iterating this step one gets that $1+1+X\in J_m=\langle X^m\rangle R$, $\forall m$. This shows that $1+1+X=0$ in $\varprojlim R_m$. Next we relate the degree of the polynomial $p(X)$ with the absolute value of the integer $p(-2)$. Let
\begin{equation}\label{bdefn}
 j(n):=\frac 13 (-2)^{n}-\frac 12
(-1)^{n}+\frac 16\qquad n\in\N.\end{equation}
The degree $n$ of a polynomial $p(X)$ with coefficients in $\{0,1\}$ specifies the sign of the integer $p(-2)$ as $(-1)^n$ and provides lower and upper bounds on the modulus $\vert p(-2)\vert$ as follows
$$
 \vert j(n-1)\vert<\vert p(-2)\vert\leq \vert j(n+1)\vert.
$$
Given an integer $m\in \Z$, the first inequality provides the following bound,  on the degree of the polynomial $p$ such that $p(-2)=m$
$$
\deg(p)\leq \log_2(3 \vert m\vert +2)+1.
$$
The projective limit $\varprojlim R_m$ is here the ring $\Z_2$ of $2$-adic integers, and the elements of $\Z$ inside $\Z_2$ are characterized by the fact that their sequence of digits is eventually constant. \vspace{.05in}

Next, we  turn to  quadratic fields for which the study of number systems in \cite{kataikovacs1,kataikovacs2} provides an exhaustive list of examples. One easily deduces from \opcit the following 
\begin{proposition}\label{kataikovacs}
    The  quadratic fields $K$ whose ring of integers admit an $\sss$-generator are 
    \begin{itemize}
    \item $\Q(\sqrt{-1})$ with generator $X=-1+\sqrt{-1}$ of the ring $\Z[\sqrt{-1}]$ of integers of $K$.
		\item $\Q(\sqrt{-2})$ with generator $X=\sqrt{-2}$ of the ring $\Z[\sqrt{-2}]$ of integers of $K$.
  \item $\Q(\sqrt{-7})$ with generator $X=\frac{1}{2} \left(1+ \sqrt{-7}\right)$ of the ring  of integers of $K$.
  \end{itemize}
\end{proposition}
\proof The norm $N(\alpha)$ of an $\sss$-generator is equal to $2$, thus the set  $N_0(\alpha):=\{0,\ldots, N(\alpha)-1\} $ defining a canonical number system in the sense of \cite{kataikovacs1,kataikovacs2} is $\{0,1\}$ and the result follows from Theorem 1 of \cite{kataikovacs2} in the complex case and Satz 1 of \cite{kataikovacs1} in the real case.\endproof

\subsubsection{The polynomial ring $(\Z[i],-1+i)$}\label{sexample0}
\begin{figure}[H]	\begin{center}
\includegraphics[scale=0.35]{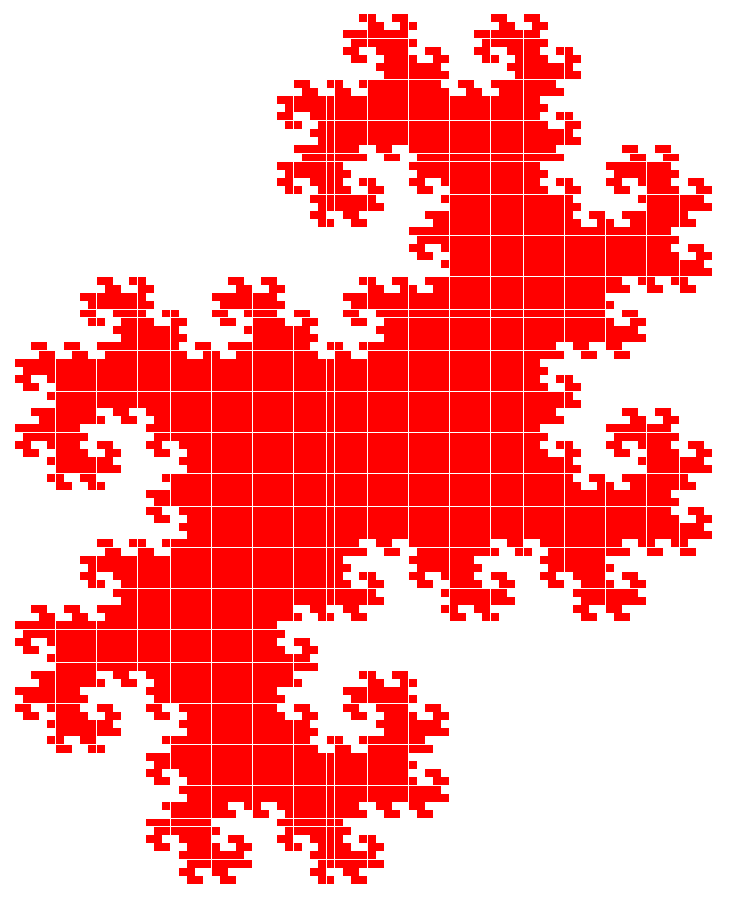}
\end{center}
\caption{Gaussian integers as $\sss$-polynomials in  degree $\leq 12$\label{gauss1}}
\end{figure}

Here, we consider the ring $R=\Z[i]$ of Gaussian integers (sometimes called binarions; see \cite{gauss}) with $X=-1+i$ as  $\sss[\mu_{1,+}]=\sss$-generator.  Indeed, every Gaussian integers can be written uniquely as a finite sum of powers of $X$ (\cite{Gilbert, Robert} and Figure \ref{gauss1}). One has the equality $1+1=P(X)=X^2+X^3$, which allows one to compute the sum of any pair of polynomials with coefficients in $\{0,1\}$.
 
\begin{proposition}\label{gcase} Let $R=\Z[i]$, $X=-1+i$. \newline
$(i)$~The ideal of $R=\Z[i]$ generated by $X^2$ is the same as the ideal generated by $2$.\newline
$(ii)$~The ring $R_m$ for $m=2k$ is $\Z/(2^k\Z)[X]$ where $X^2=-2-2X$.\newline
$(iii)$~The ring $R_m$ for $m=2k+1$ is $\Z/(2^{k+1}\Z)\oplus \Z/(2^{k}\Z)\, X$ where $X^2=-2-2X$.\newline
$(iv)$~The projective limit  $\varprojlim R_m$ is the ring $\Z_2[i]\sim \Z_2[X]$  where $X^2=-2-2X$.    
\end{proposition}
\proof $(i)$~The element $U=(1+X)\in R$ is a unit since $U^4=1$ and one has
$$
X^2=-2-2X\in 2R, \ \ 2=-(1+X)^{-1}X^2\in X^2R
$$
$(ii)$~By $(i)$ the ideal $X^{2k}R$ is equal to $2^k R$. One has $R=\Z[i]$ and $R/(2^kR)=\Z/(2^k\Z)[i]=\Z/(2^k\Z)[X]$ with  $X^2=-2-2X$, thus one gets $(ii)$. \newline
$(iii)$~Let $m=2k+1$. Any element of $R$ is of the form $z=a+ b X$ where $a,b\in \Z$. In $R$ one has $2^{k+1}\in X^{2k+2}R\subset J_m$ and $2^k X\in X^{2k+1}R= J_m$. Thus the homomorphism $\Z[X]\to R_m$ induces a surjective homomorphism from $\Z/(2^{k+1}\Z)\oplus \Z/(2^{k}\Z)\, X$ to $R_m$. It is bijective since the cardinalities are equal.\newline
$(iv)$~The extension $\Q_2[i]$ is totally ramified of index $e=2$ (see \cite{Robert}, 4.2). The polynomial $X^2+2X +2$ is an Eisenstein polynomial  which defines $\Q_2[i]$ as its splitting field. The valuation of $X$ is one half of the valuation of $2$.\endproof 
\subsubsection{The polynomial ring $\left(\Z[\sqrt{-2}],\sqrt{-2}\right)$}
~The element $X:=\sqrt{-2}$ is an $\sss[\mu_{1,+}]=\sss$-generator of the  ring of integers $\Z[\sqrt{-2}]$  of the imaginary quadratic field $\Q(\sqrt{-2})$. This follows directly from \S \ref{Z2} and Proposition \ref{sprop1}. The hold is given by the polynomial $P(X)=X^4+X^2$. A straightforward analogue of Proposition \ref{gcase} holds.
\subsubsection{The polynomial ring $\left(\cO(\Q(\sqrt{-7})),\frac{1}{2} (1+ \sqrt{-7})\right)$}

\begin{figure}[H]	\begin{center}
\includegraphics[scale=0.45]{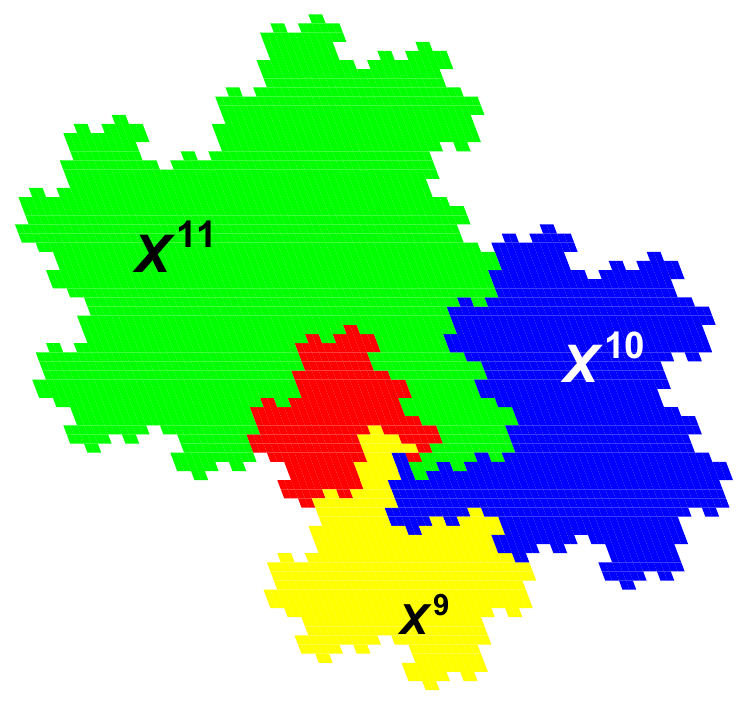}
\end{center}
\caption{Polynomials of degree $\leq 11$ for $X=\frac{1}{2} \left(1+ \sqrt{-7}\right)$\label{sevenseven}}
\end{figure}

The element $X:=\frac{1}{2} \left(1+ \sqrt{-7}\right)$ is an $\sss[\mu_{1,+}]=\sss$-generator of the  ring $\cO(\Q(\sqrt{-7}))$ of integers of the imaginary quadratic field $\Q(\sqrt{-7})$. The hold is given by the polynomial $P(X)=X^3+X$. Let $F$ be the fundamental domain of $\cO(\Q(\sqrt{-7}))$ given by the parallelogram with vertices $0,1,X,X+1$.
Figure \ref{sevenseven} shows the neighborhood of $0\in \C$ obtained as the union of the  translations $F+p(X)$ by polynomials $p(X)$ of degree $\leq 11$. 
\begin{proposition}\label{7case} Let $R=\cO(\Q(\sqrt{-7}))$, $X=\frac{1}{2} \left(1+ \sqrt{-7}\right)$. \newline
$(i)$~The ring $R_m$  is $\Z/(2^m\Z)$. \newline
$(ii)$~The projective limit  $\varprojlim R_m$ is the ring $\Z_2$.  \newline
$(iii)$~The element $X\in\varprojlim R_m=\Z_2$ is the only solution divisible by $2$ in the ring $\Z_2$ for the equation $2+X+X^2=0$.
\end{proposition}
\proof The hold is given by $P(X)=X^3+X$ and one has $P(X)-2=(X-1) \left(X^2+X+2\right)$. By Hensel's Lemma, the equation $2+X+X^2=0$ admits a unique solution $\alpha$ in $\Z_2$ of the form $\alpha=1+2\epsilon$ and a unique solution of the form  $\beta=2(1+2\epsilon')$. In fact one has $\alpha \beta=2$ and $\alpha+\beta=-1$.  The homomorphism $\rho:\Z[\frac{1}{2} \left(1+ \sqrt{-7}\right)]\to \Z_2$ given by $\rho\left(\frac{1}{2} \left(1+ \sqrt{-7}\right)\right)=\beta$ is well defined since $\beta$ is a solution of the equation $2+X+X^2=0$. Moreover $\beta$ is the product of $2$ by a unit of $\Z_2$ (but this fails in $R=\cO(\Q(\sqrt{-7}))$). The projection $X_m$ of $\beta$ in $\Z_2/(2^m\Z_2)=\Z/(2^m\Z)$ fulfills $P(X_m)=2$ and $X_m$ is the product of $2$ by a unit. Thus the ideals generated by powers of $X_m$ are the same as those generated by powers of $2$. This proves the three assertions $(i)$, $(ii)$, $(iii)$. \endproof

\subsection{Polynomial rings in one generator over  $\spm$}
\subsubsection{The polynomial ring $(\Z,3)$}~The case of the $\spm$-generator $3\in \Z$ is particularly relevant because, as shown in \cite{RR}, the addition coincides with that of the Witt vectors in $W(\F_3)=\Z_3$.
\begin{proposition}\label{zcase} Let $R=\Z$, $X=3$ is an $\spm$-generator of $R$. The hold is $P(X)=-1+X$. \newline
$(i)$~The ring $R_m$  is $\Z/(3^m\Z)$. \newline
$(ii)$~The projective limit  $\varprojlim R_m$ is the ring $W(\F_3)=\Z_3$.  \newline
$(iii)$~The set of Witt vectors with only finitely many non-zero components forms a subring of $W(\F_3)$ isomorphic to $\Z$.\end{proposition}
In order to organize the next examples we give the  list of imaginary  quadratic field extensions of $\Q$ generated by rings of $\spm$-polynomials in one variable.
\begin{proposition}\label{ringquad} 	The imaginary quadratic fields $K$ generated by rings of $\spm$-polynomials in one variable are 
\begin{itemize}
		\item $\Q(\sqrt{-2})$ with generator $X=1+\sqrt{-2}$ of the ring $\Z[\sqrt{-2}]$ of integers of $K$.
	\item $\Q(\sqrt{-3})$ with generator $X=\sqrt{-3}$ of $\Z[\sqrt{-3}]$ (not a UFD).
	\item $\Q(\sqrt{-11})$ with generator $X=\frac 12(1+\sqrt{-11})$ of the ring  of integers of $K$.
\end{itemize}	
\end{proposition}
\proof Let $P(X)=-1+\sum_{j=1}^{n-1} a(j)X^j+\epsilon X^n$, $\epsilon \in \{\pm 1\}$, $a(j) \in \{-1,0,1\}$,  be the carry leading to an imaginary quadratic extension. The roots of the polynomial $P(X)-2$ are algebraic integers, and we assume that one of them, say $\alpha$, is quadratic imaginary. Let $q(x)=x^2-b x+c$ be its minimal polynomial. It has integral coefficients so $b,c\in \Z$, and by definition, it divides $P(X)-2$. The constant coefficient $c$ must be equal to $3$. Indeed it divides the constant coefficient $-3$ of $P(X)-2$ and since $b^2-4c<0$ it is positive. It cannot be equal to $1$ since in that case one would get  $b\in \{-1,0,1\}$, and $\alpha\in \{i,j,-j\}$ which contradicts the injectivity of the map $\sigma$. For $c=3$ the possible values of $b$ are $b=0$ which gives the solution $\alpha=\sqrt{-3}$, $b=\pm 1$ which gives the solutions $\alpha= \frac{1}{2}\left(\pm 1\pm i \sqrt{11}\right)$,  $b=\pm 2$ which gives the solutions $\alpha= \pm 1\pm i \sqrt{2}$,  and finally $b=\pm 3$. We shall now show that this last choice which gives $\alpha =\frac{1}{2} \left(\pm 3\pm i \sqrt{3}\right)$ does not give a solution. To prove this it is enough to show that the polynomial $3+3X +X^2$ cannot divide a polynomial $P(X)-2$ with $P$ of the above form. We thus assume an equality of the form 
$$
(3+3X +X^2)\left(\sum_{j=0}^{n-2} b(j)X^j\right)=-3+\sum_{j=1}^{n-1} a(j)X^j+\epsilon X^n,  \ \epsilon \in \{\pm 1\}, \ a(j) \in \{-1,0,1\}
$$
Since the coefficients of $P-2$ are integers and the leading coefficient of $3+3X +X^2$ is $1$ the coefficients $b(j)$ are integers.  We get $b(0)=-1$, $3 b(1)-3=a(1)$, but $a(1)\in \{-1,0,1\}$ and thus working modulo $3$ one gets $a(1)=0$ and hence $b(1)=1$. Considering the coefficient of $X^2$ we get $3 b(1)+3 b(2)-1=a(2)$ which gives $a(2)=-1$ and $b(2)=-b(1)=-1$. We can now work by induction to show that $b(j)=(-1)^{j+1}$. Indeed the coefficient of $X^j$ is 
$b(j-2)+3 b(j-1)+3 b(j)=a(j)$ and if we know that $b(j-2)=(-1)^{j-1}$ and $b(j-1)=(-1)^{j}$ we get $a(j)=b(j-2)$ and $3 b(j-1)+3 b(j)=0$ so that $b(j)=(-1)^{j+1}$. This works for $j\leq n-2$. The coefficient of $X^{n-1}$ is $b(n-3)+3b(n-2)=a(n-1)$ and this gives a contradiction since one gets $a(n-1)=b(n-3)$ (working modulo 3) which contradicts the fact that $b(n-2)\neq 0$. \endproof
\subsubsection{The polynomial ring $(\cO(\Q[\sqrt{-11}]),\frac{1}{2} \left(1+ \sqrt{-11}\right))$}~This section is dedicated to a detailed proof that $X:=\frac{1}{2} \left(1+ \sqrt{-11}\right)$ is an $\spm$-generator of the ring of integers of the number field $\Q(\sqrt{-11})$. The reason for providing the details of the proof is because we want  to emphasize that in such a case, and unlike working over $\sss$, one can explicitly control the cancellations in the computations. 
\begin{proposition}\label{ring11} Let $\cO$ be the ring of integers of the number field $\Q(\sqrt{-11})$. 
\newline $(i)$~$X:=\frac{1}{2} \left(1+ \sqrt{-11}\right)$ is an $\spm$-generator of $\cO$. The hold of $(\cO,X)$ is $P(X)=-1+X-X^2$.	\newline
$(ii)$~The projective limit  $\varprojlim R_m$ is the ring $W(\F_3)=\Z_3$.
\end{proposition}
 \begin{figure}[H]	\begin{center}
\includegraphics[scale=0.4]{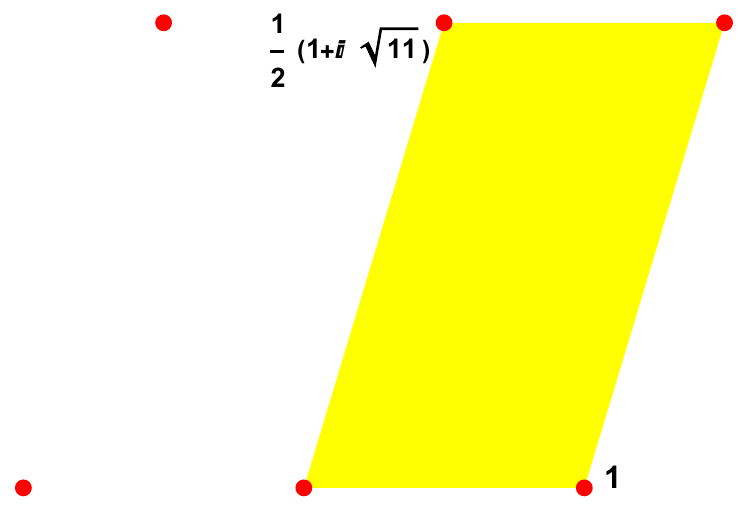}
\end{center}
\caption{Fundamental domain of the lattice $\cO$\label{fund-domain}}
\end{figure}
\no The proof  requires  a preliminary lemma. We first recall some classical results concerning the ring of integers $\cO$ of the imaginary quadratic field $K=\Q(\sqrt{-11})$. The discriminant of $K$ is $d=-11$. Thus since $-11\sim 1$ modulo $4$, the lattice $\cO$ is  $\Z+\Z X$ where $X:=\frac{1}{2} \left(1+\sqrt{-11}\right)$. By construction one has
\begin{equation}\label{basic}
1+1=P(X), \qquad P(X)=-1+X-X^2.
\end{equation}
One wants to show that every element  $z\in \cO$ can be  written uniquely as a polynomial $z=\sum_j \alpha_j\, X^j$, with $\alpha_j\in \{-1,0,1\}$. Figure \ref{fund-domain}   shows the translates of the fundamental domain of the lattice, while the next figures provide a sketch of a few steps of the process of   representing elements of $\cO$ in terms of polynomials of degree $\leq n$,  showing those described by polynomials of degree $=n$ with a new color.
\begin{figure}[H]
\centering
\begin{subfigure}{.5\textwidth}
  \centering
\includegraphics[width=.7\linewidth]{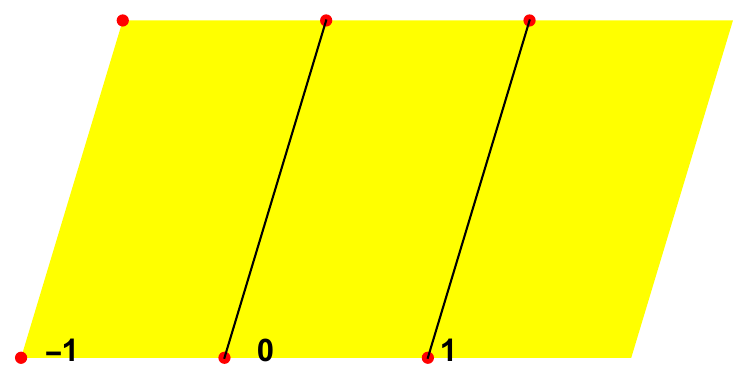}
  \caption{First step, polynomials of degree $0$}
  \label{fund1}
\end{subfigure}%
\begin{subfigure}{.5\textwidth}
  \centering
\includegraphics[width=.45\linewidth]{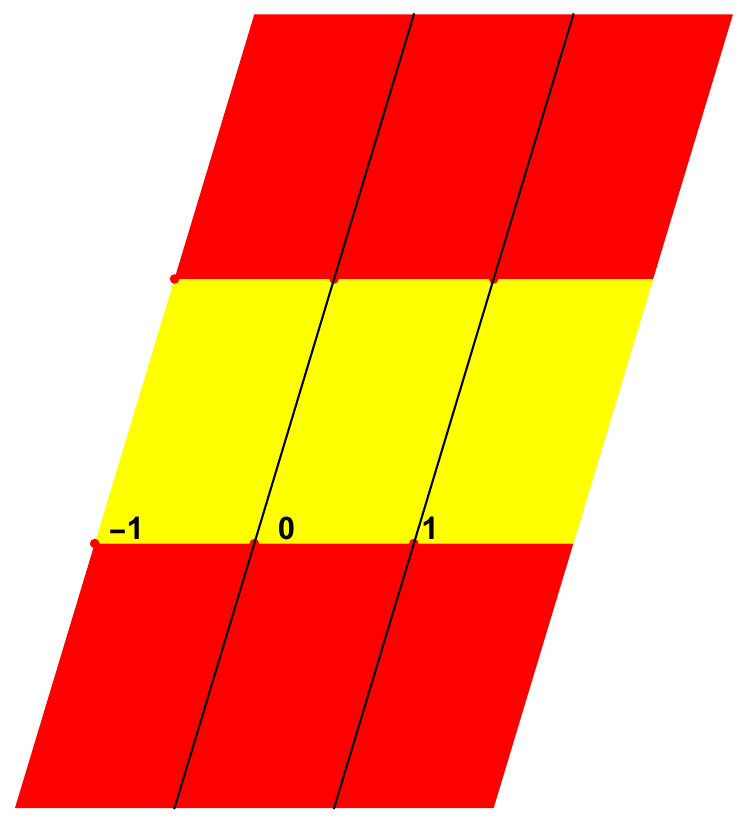}
  \caption{Second step, polynomials of degree $\leq 1$}
  \label{fund2}
\end{subfigure}
\caption{The first two steps}
\end{figure}
\begin{figure}[H]
\centering
\begin{subfigure}{.45\textwidth}
  \centering
\includegraphics[width=.7\linewidth]{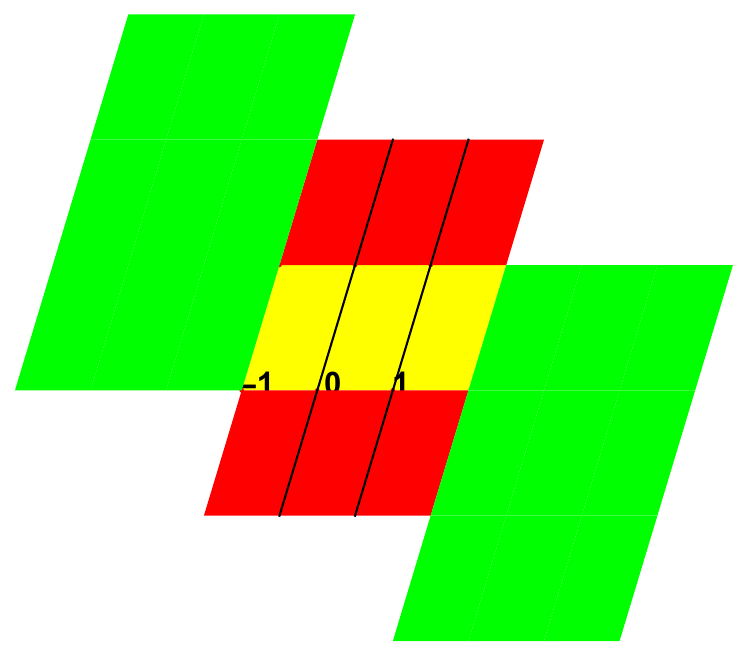}
  \caption{Third step, polynomials of degree $\leq 2$}
  \label{fund3}
\end{subfigure}%
\begin{subfigure}{.5\textwidth}
  \centering
\includegraphics[width=.7\linewidth]{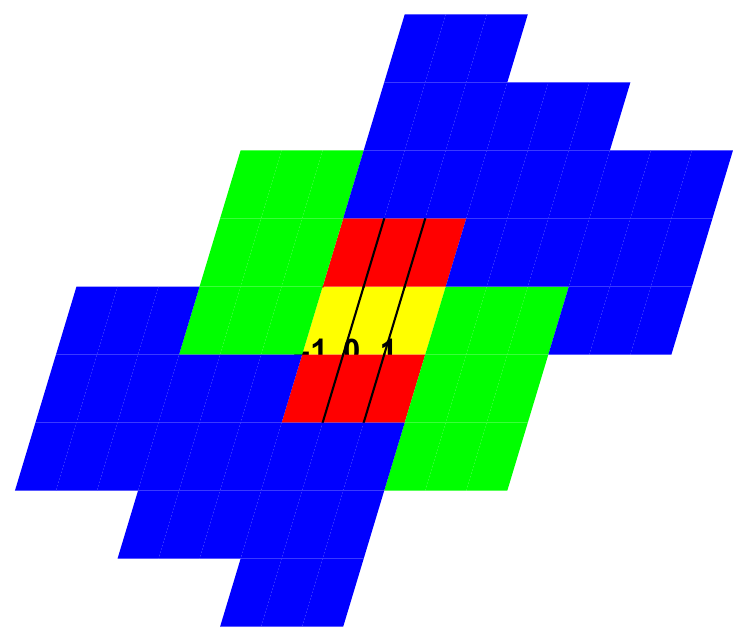}
  \caption{Fourth step, polynomials of degree $\leq 3$}
  \label{fund4}
\end{subfigure}
\caption{The third and fourth  steps}
\end{figure}
\begin{figure}[H]
\centering
\begin{subfigure}{.45\textwidth}
  \centering
\includegraphics[width=.6\linewidth]{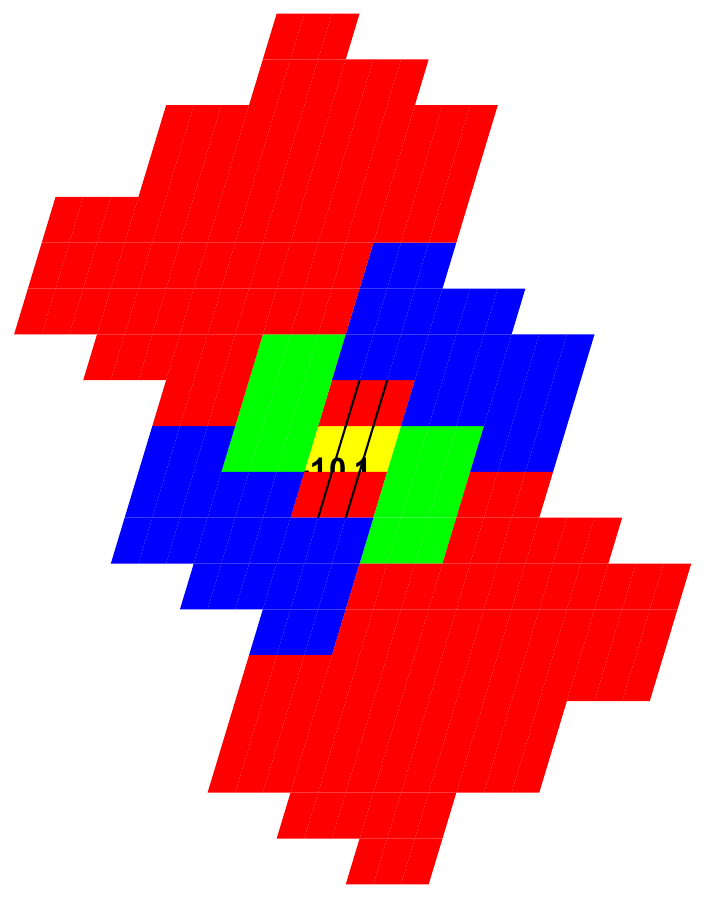}
  \caption{Fifth step, polynomials of degree $\leq 4$}
  \label{fund5}
\end{subfigure}%
\begin{subfigure}{.45\textwidth}
  \centering
\includegraphics[width=.85\linewidth]{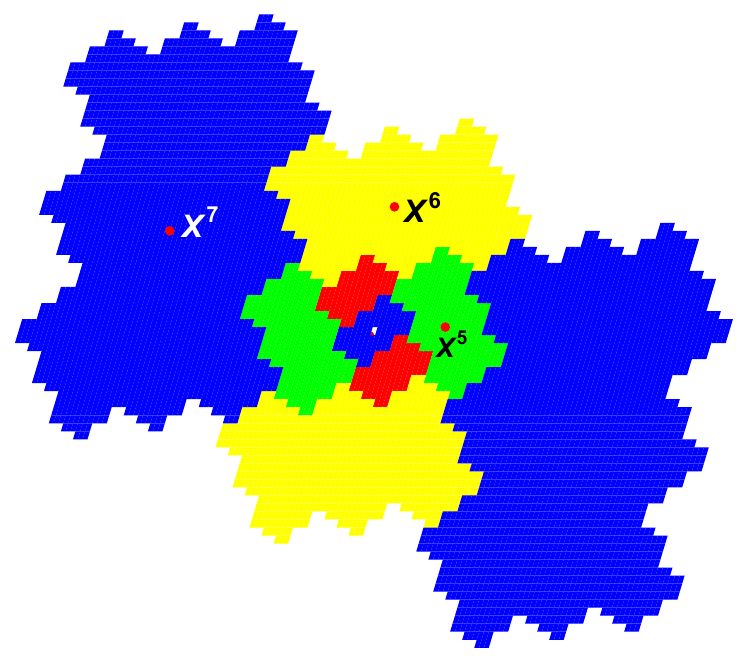}
  \caption{Eigth's step, polynomials of degree $\leq 7$}
  \label{fund6}
\end{subfigure}
\caption{The fifth and eighth steps}
\end{figure}

\no By comparing Figures \ref{fund1}, \ref{fund2}, \ref{fund3}, \ref{fund4}, \ref{fund5}, \ref{fund6}, one notices that the translation $z\mapsto z+1$ does not increase the degree of the polynomial by more than $2$ units. Next lemma provides a formal proof of this fact.
\begin{lemma}\label{atmost2} Let $z=\sum_{j=0}^n \alpha_j X^j\in\cO$, $\alpha_j\in \{-1,0,1\}$. Then there exist coefficients $\beta_j\in \{-1,0,1\}$, with $0\leq j\leq n+2$, such that $z+1=\sum_{j=0}^{n+2} \beta_j X^j$.	
\end{lemma}
\proof We proceed by induction on the integer $n$. For $n=0$, the result follows from \eqref{basic}. Let us assume that the result is proved up to $n-1$, then there exists coefficients $\gamma_j \in \{-1,0,1\}$ such that
$$
z=\left(\sum_{j=0}^{n-1} \alpha_j X^j\right)+\alpha_nX^n ~\Longrightarrow~ z+1=
\left(\sum_{j=0}^{n+1} \gamma_j X^j\right)+\alpha_nX^n.
$$
Let us consider a sum such as $\gamma_nX^n+\gamma_{n+1}X^{n+1}+\alpha_nX^n$ and express it without going beyond $X^{n+2}$. If $\gamma_{n+1}=0$ this follows again from \eqref{basic}. We can thus assume that $\gamma_{n+1}=\pm 1$ and also that both $\gamma_n$ and $\alpha_n$ are non-zero and equal since otherwise the sum $\gamma_nX^n +\alpha_nX^n$ would have  degree at most $n$. The only case to exclude then is when  $\gamma_n$, $\alpha_n$, and $\gamma_{n+1}$ are all equal (and non-zero), since only in that case would one get a term in $X^{n+3}$ from the sum 
\begin{align*}
&X^n+X^n+X^{n+1}=X^n(1+1+X)=X^n(-1+X+X-X^2)=\\
&=X^n(-1-X+X^2-X^2-X^3)=-X^n-X^{n+1}-X^{n+3}.
\end{align*}
To exclude this case, one adds to the induction hypothesis the condition that if the last term $\beta_{n+2}$ of the polynomial of degree $n+2$ representing $z+1$ is non-zero, then the term $\beta_{n+1}$ is zero or of the opposite sign. This condition is fulfilled for $n=0$, and if we assume it for $n-1$, it holds also for $n$. Indeed, the only cases when $\beta_{n+2}\neq 0$ 
arise when either $\gamma_{n+1}=0$, in which case $\beta_{n+1}$ and $\beta_{n+2}$ have opposite signs, or $\gamma_{n+1}=\epsilon=\pm 1$ in which case $\gamma_n=\alpha_n=-\epsilon$, which gives 
$$
\gamma_nX^n+\gamma_{n+1}X^{n+1}+\alpha_nX^n=-\epsilon X^n(1+1-X)=\epsilon X^n+\epsilon X^{n+2},
$$
implying that $\beta_{n+1}=0$ in this case. Thus the induction hypothesis still holds for $n$, and this concludes the proof. \endproof

\proof(of Proposition \ref{ring11})~Lemma \ref{atmost2}  holds for the abstract law of addition defined using \eqref{basic} on the projective limit of the $R_n$. The proof  shows that elements of this limit, which have only a finite number of non-zero coordinates, are stable under the addition of $1$. Using \eqref{polX1}, it follows that they are also stable under the addition of any monomial and hence that they form an additive group $A$. Thus, it remains to show that the map $\rho:A\to \C$ defined by
$$
\rho\Big(\sum_j \alpha_j X^j\Big):=\sum_j \alpha_j z^j, \qquad z=\frac{1}{2} \left(1+\sqrt{-11}\right)
$$ 
is injective. Let $\sum_j \alpha_j X^j\in \ker \rho$, then  $\sum_j \alpha_j z^j=0$ and thus $z$ fulfills an equation $E(z)=0$ with integral coefficients whose leading coefficient is $1$ and the constant term is $\pm 1$. The polynomial $E$ is thus a multiple of the minimal polynomial $z^2-z+3$ of the field extension. The quotient polynomial has integral coefficients; thus, one gets a contradiction using the product of constant terms.\endproof
\subsubsection{The polynomial ring $\left(\Z[\sqrt{-3}],\sqrt{-3}\right)$}
~The element $X:=\sqrt{-3}$ is an $\sss[\mu_{2,+}]=\spm$-generator of the  ring $\Z[\sqrt{-3}]$ and the latter is a maximal order in the ring  of integers of the imaginary quadratic field $\Q(\sqrt{-3})$. This follows directly from \S \ref{Z2} and Proposition \ref{sprop1}. The hold is given by the polynomial $P(X)=-1-X^2$. A straightforward analogue of Proposition \ref{gcase} holds.
\subsubsection{The polynomial ring $\left(\cO(\Q(\sqrt{-2})),1+ \sqrt{-2}\right)$}~One obtains  similarly that $P(X)=-1-X+X^2-X^3$ is the hold associated to the $\spm$ generator $1+ \sqrt{-2}$ of the ring of integers of the imaginary quadratic field $\Q(\sqrt{-2})$
\begin{figure}[H]	\begin{center}
\includegraphics[scale=0.4]{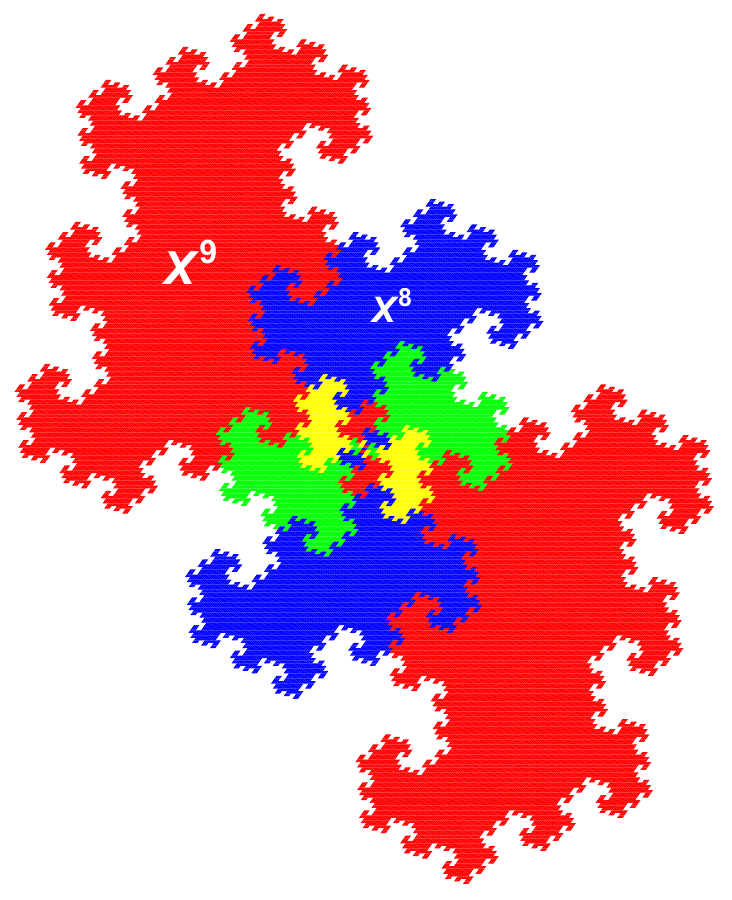}
\end{center}
\caption{Polynomials of degree $\leq 9$ for $X=1+i \sqrt{2}$\label{rootofmin2}}
\end{figure}
\begin{proposition}\label{ringi2} Let $\cO$ be the ring of integers of the number field $\Q(\sqrt{-2})$. 
\newline $(i)$~$X:=1+ \sqrt{-2}$ is an $\spm$-generator of $\cO$. The hold of $(\cO,X)$ is $P(X)=-1-X+X^2-X^3$.	\newline
$(ii)$~The projective limit  $\varprojlim R_m$ is the ring $W(\F_3)=\Z_3$.
\end{proposition}
\no  Figure \ref{rootofmin2} reproduces the pattern obtained by inputting polynomials of degree $\leq 9$. In this case, the analog of Lemma \ref{atmost2} holds with the bound $n+3$ instead of $n+2$.


\subsection{Polynomial rings in one generator over  $\sss[\mu_{3,+}]$}
In the next example  the field $R_1$ is the finite field $\F_4$. One lets $\mu_{3,+}\subset \C$ be the solutions of $x(x^3-1)=0$, $j=\exp(2\pi i/3)$ and $\Z(j)\subset \Q(j)$ be the ring of integers of the quadratic imaginary field $\Q(j)$.
\begin{proposition} \label{mu3} $(i)$~The number $-2\in \Z(j)$ is an $\sss[\mu_{3,+}]$-generator of the ring $R=\Z(j)$.\newline
$(ii)$~The hold is given by 
$$
h(1)=X+X^2, \ \ h(j)=j^2 X+j^2, \ \ h(j^2)=j X+j
$$
$(iii)$~The  field $R_1$ is the finite field $\F_4$.\newline
$(iv)$~The projective limit  $\varprojlim R_m$ is the Witt ring $W(\F_4)$ and the ring $R_m$ is the quotient of $W(\F_4)$ by $2^m\,W(\F_4)$.    
\end{proposition}

\proof Let $J=2\Z(j)\subset \Z(j)$, then $J^n$ is the ideal generated by $X^n$ where $X=-2$. Let $\sigma:\cP(\mu_3)\to R=\Z(j)$ be the map defined by \eqref{sigma}. For each $n$ the composition $\pi_n\circ \sigma$, from the subset  $\cP^{n-1}(\mu_3)\subset\cP(\mu_3)$ formed of polynomials of degree $<n$ to the quotient ring $R_n=R/J^n$, is surjective and hence injective since the cardinalities of source and target are the same. It follows that the map $\sigma:\cP(\mu_3)\to R=\Z(j)$ is injective. To show that it is surjective one uses the general method involving the limit of the subsets 
$$
Z_n:=(-2)^{-n}\left(\sigma(\cP^{n}(\mu_3)+F)\right)\subset \C
$$
where $F$ is a fundamental domain for $\Z(j)$. One observes that passing from $n$ to $n+1$ only alters $Z_n$ on its boundary and that $Z_n$ contains an open disk centered at $0$. 

\begin{figure}[H]	\begin{center}
\includegraphics[scale=0.35]{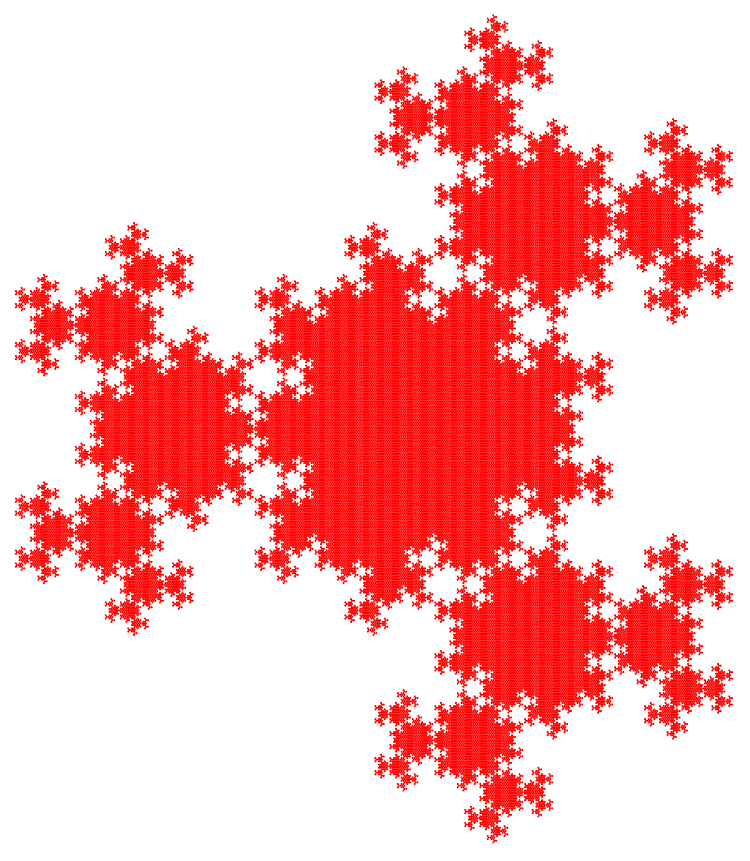}
\end{center}
\caption{Polynomials of degree $\leq 7$ for $X=-2$\label{smu3}}
\end{figure}
\subsection{Polynomial rings in one generator over  $\sss[\mu_{4,+}]$}
\begin{figure}[H]	\begin{center}
\includegraphics[scale=0.4]{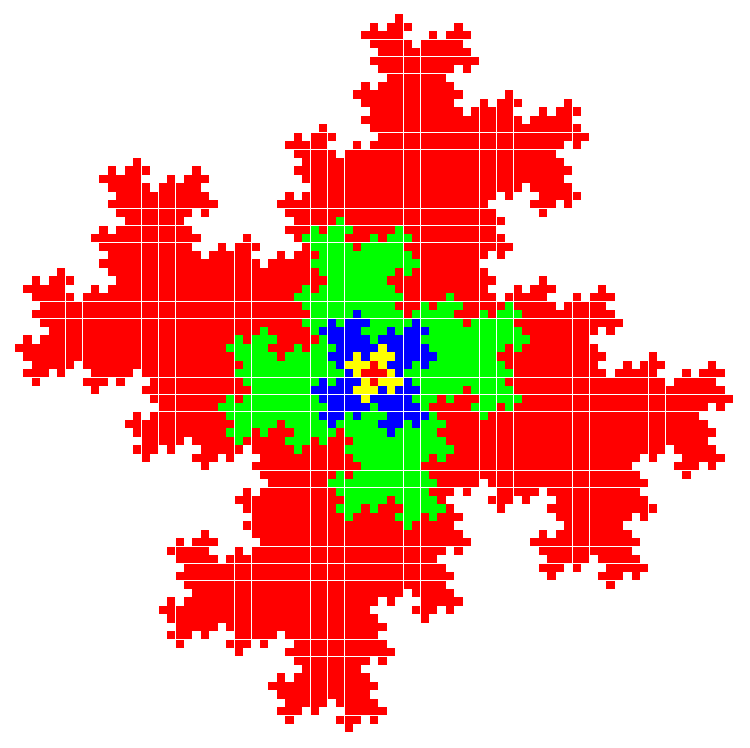}
\end{center}
\caption{Polynomials of degree $\leq 4$ for $X=1+2i$\label{gauss4}}
\end{figure}

\begin{proposition} \label{mu4} $(i)$~The number $X=1+2i$ is an $\sss[\mu_{4,+}]$-generator of the ring $R=\Z(i)$.\newline
$(ii)$~The hold is given by  $h(0)=1$ and
$$
h(1)=i - i\, X, \ \ h(i)=-i+ X , \ \ h(-i)=-1-i\, X
$$
$(iii)$~The  field $R_1$ is the finite field $\F_5$.\newline
$(iv)$~The projective limit  $\varprojlim R_m$ is the Witt ring $W(\F_5)=\Z_5$ and the ring $R_m$ is the quotient of $W(\F_5)$ by $5^m\,W(\F_5)$.    
\end{proposition}
\proof In the $p$-adic field $\Z_5$ there exists a unique square root of $-1$ equal to $2$ modulo $5$ (see \cite{Robert}, \S 6.7). Let $\rho:\Z(i)\to \Z_5$ be the unique morphism such that, modulo $5$, one has $\rho(i)=2$. Then $\rho(X)=5 u$ where $u$ is a unit in $\Z_5$. The morphism $\rho$ restricted to  
$\mu_{4,+}=\{0,1,i,-1,-i\}$ gives a multiplicative section of the quotient map $R\to R/XR$. One has $\Z_5/\rho(X)^m\Z_5=\Z/5^m\Z$ and the morphism $\rho$ induces an isomorphism $R_m\simeq \Z_5=\Z/5^m\Z$. Statements $(iii)$ and $(iv)$ follow, as well as the injectivity of 
 the map 
 $\sigma:\cP(\mu_4)\to R=\Z(i)$. One can prove the surjectivity of $\sigma$ as for  
Proposition \ref{mu3} using Figure \ref{gauss4}. Statements $(i)$ and $(ii)$ follow. \endproof

\subsection{Polynomial rings in one generator over  $\sss[\mu_{6,+}]$}

\begin{proposition} \label{muj6} $(i)$~The number $X=2-j$ is an $\sss[\mu_{6,+}]$-generator of the ring $R=\Z(j)$.\newline
$(ii)$~The hold is given by $h(j)=j+1$, $h(j^2)=j^2+1$, $h(0)=1$ and
$$
h(1)=X+j, \ \ h(-j^2)=-j^2\, X+j^2, \ \ h(-j)=-1+X
$$
$(iii)$~The  field $R_1$ is the finite field $\F_7$.\newline
$(iv)$~The projective limit  $\varprojlim R_m$ is the Witt ring $W(\F_7)=\Z_7$ and the ring $R_m$ is the quotient of $W(\F_7)$ by $7^m\,W(\F_7)$.    
\end{proposition}
The proof can be easily deduced from \cite{Robert}, \S 4.6.  

\begin{figure}[H]	\begin{center}
\includegraphics[scale=0.4]{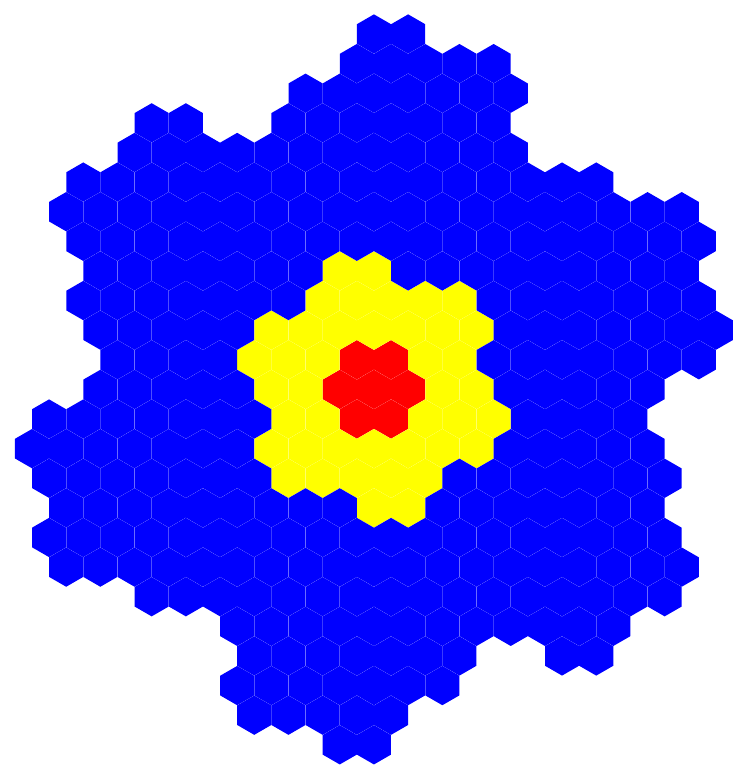}
\end{center}
\caption{Polynomials of degree $\leq 2$ for $X=\frac{5}{2}-\frac{i \sqrt{3}}{2}$\label{hexa6}}
\end{figure}

\vspace{.2in}

\begin{minipage}{0.60\linewidth}
\text{Alain Connes}\par
\textsc{Coll\`ege de France}\par
3 Rue d'Ulm\par F-75005 Paris, France\par
IHES\par
35 Rte de Chartres\par 
91440 Bures-sur-Yvette, France\par
Email: \texttt{alain@connes.org}
\end{minipage}
\hfill
\begin{minipage}{0.6\linewidth}
\text{Caterina Consani}\par
Department of Mathematics\par
\textsc{The Johns Hopkins University}\par
3400 N Charles Street\par
Baltimore MD 21218, USA\par
Email: \texttt{cconsan1@jhu.edu}\par
~\par
~
\end{minipage}
\end{document}